\documentclass[10pt]{article}

\usepackage{bbm}
\usepackage{bbold}
\usepackage{amsmath}
\usepackage{amsthm}
\usepackage{amssymb} 
\usepackage{mathrsfs}
\usepackage{graphicx}
\usepackage{authblk}
\usepackage{framed}
\usepackage[backref=none,hypertexnames=false,colorlinks=true,urlcolor=blue,linkcolor=blue,citecolor=blue]{hyperref}
\usepackage[letterpaper,text={7in,9in},centering]{geometry}
\usepackage{comment}
\usepackage[table,xcdraw]{xcolor}
\usepackage[utf8]{inputenc}
\usepackage{graphicx}
\usepackage{float}
\usepackage{caption}
\usepackage{subcaption}
\usepackage{endnotes}
\usepackage{matlab-prettifier}

\newcommand{\red}[1]{\textcolor{black}{#1}}

\usepackage[mathlines]{lineno}

\usepackage[english]{babel}


\def\rbar{\overline{r}}
\def\gbar{\overline{g}}

\def\ri{{\rm i}}

\def\d{{\rm d}}

\def\rb{{\bf r}}
\def\bX{{\bf X}}
\def\Vb{{\bf V}}
\def\that{\widehat{t}}
\def\Phihat{\widehat{\Phi}}
\def\Psihat{\widehat{\Psi}}
\def\Ghat{\widehat{G}}
\def\what{\widehat{w}}
\def\Vhat{\widehat{V}}

\def\Re{{\rm Re}}

\usepackage{accents}
\newcommand*{\bigdot}[1]{%
  \accentset{\mbox{\large\bfseries .}}{#1}}

\begin{document}

\title{Dynamic interactions and equilibrium configurations of pulses in the two-dimensional complex quintic Ginzburg-Landau equation}
\author{M.R.Turner\footnote{Corresponding author: m.turner@surrey.ac.uk}\,  and  D.J.B.Lloyd}

\affil{\small School of Mathematics and Physics, University of Surrey, Guildford, GU2 7XH, UK}

\date{}
\maketitle






\begin{abstract}
This paper constructs a fast and effective novel numerical scheme which accurately calculates the dynamics of weakly-interacting pulses in the two-dimensional quintic-complex Ginzburg-Landau equation (QCGLE). The numerical scheme uses a global centre-manifold reduction, where the solution to the QCGLE is constructed as the sum of the individual pulses plus a remainder function, which is chosen to be orthogonal to the zero adjoint eigenmodes of the QCGLE linear operator. Projecting this constructed solution onto the stable centre-manifold leads to a fast-slow system of equations consisting of {\it slow} ordinary differential equations for the position and phases of the individual pulses and a {\it fast} partial differential equation for the remainder function. By considering the pulses to be well-separated, the system can be expanded asymptotically in terms of the small parameter $\epsilon=e^{-\lambda_r d}\ll1$, where $\lambda_r$ is the spatial decay rate of the pulse, and $d>0$ is the minimum pulse separation distance. Here the remainder function is determined via a stationary partial differential equation that can be readily solved in an efficient manner using GMRES. Results for $N=2,3,4$ and $5$ pulses are considered, and it is found that different equilibrium solutions are possible such as stable fixed points and limit cycles. The interaction of two stable $N=3$ coherent structures is also considered, where the common tendency found is for the structure to degenerate into pairs of pulses which propagate away from the initial configuration of pulses.
\end{abstract}

\section{Introduction}

Interest in the existence and interaction of spatially localized structures in physical systems is vast; see for instance the various reviews~\cite{knobloch2015spatial,knobloch2008spatially,Malomed2022,bramburger2025localizedpatterns}. One of the most common localized structures to be observed in these systems is the soliton or pulse, also known as a spot. Such structures have a central peak surrounded by exponentially decaying tails back to the background state, and have been observed in a variety of systems such as vibrating granular materials \cite{umbanhowar1996}, shear flows \cite{duguet2009,schneider2010}, neural networks \cite{laing2002,laing2003}, crime hot-spots in criminological models \cite{lloyd2013}, solitons in nonlinear optics \cite{vladimirov2002,ackemann2009} and ferrofluids \cite{richter2005,lloyd2015}. These localized pulse solutions occur in systems due to balances between nonlinear effects and either, dissipation or dispersion \cite{ablowitz1991,aranson2002}. These structures have also been extensively studied mathematically using canonical differential equations such as the Swift-Hohenberg equation \cite{burke2006,lloyd2009}, the Ginzburg-Landau equation \cite{sakaguchi2001,ZelMiel09,rossides2023}, the nonlinear Schr{\"o}dinger equation \cite{rogers2012,choudhuri2016}, the Klein-Gordon equation \cite{perel2003} and the Kadomtsev-Petviashvili equation \cite{whang2016} to name a few. When these pulses are sufficiently separated, then their dynamics are solely due to their individual construction, however as the distance between the pulses reduces they begin to interact via their exponential tails forming so-called {\it coherent structures} or {\it clusters} \cite{tlidi1994,vladimirov2001,vladimirov2002,Turaev07}. \red{If the pulses become closer still, then the interacting dynamics become strong, and the pulses can behave like elastic objects or even collide. A collision between pulses can lead to annihilation, merging into a single spot, or even more general complex dynamics such as spatiotemporal chaos \cite{nishiura2001,nishiura2003,nishiura2003chaos,nishiura2005}.} The main focus of this paper is to investigate the possible dynamics of \red{weakly interacting pulses} in the two-dimensional quintic complex Ginzburg-Landau equation.

As the magnitude of the nonlinear interaction between pulses is exponentially small with respect to their separation distance,  numerically calculating their interaction dynamics is highly challenging. Direct numerical methods are usually ineffective due to the fact that the interactions are small and the time-scales are typically long, leaving the integration scheme prone to numerical round-off errors. A more effective numerical scheme to investigate these delicate dynamics was first investigated by \cite{rossides2023} and was based on the so-called {\it centre-manifold reduction} approach \cite{ei2002,Bjorn02,ZelMiel09}. The approach is valid for any dissipative partial differential equation of the form
\begin{equation}
\frac{\partial u}{\partial t}=\mathbb{A}u+f(u),~~~~~f(0)=f'(0)=0,
\label{eqn:disppde}
\end{equation}
where $u({\bf x},t)$ is an unknown vector function, $\mathbb{A}$ is an elliptic differential operator in ${\bf x}\in\mathbb{R}^m$, $t\in[0,\infty)$ is time and $f$ is a given smooth, nonlinear function. If it can then be shown that (\ref{eqn:disppde}) exhibits an equilibrium pulse  solution $V({\bf x})$ then it in fact possesses a whole manifold of solutions, parameterised by $\boldsymbol\xi(t)$, due to system symmetries, where $\boldsymbol\xi(t)$ lives in the Lie group of the system symmetries. A solution is then sought in terms of a sum of these individual pulses, plus a remainder function $w({\bf x},t)$, which is assumed to be small compared to the pulse magnitudes. It is then shown that the dynamics of the system can be written as a fast-slow system; A set of {\it slow} ordinary differential equations (ODEs) for the Lie group parameters $\boldsymbol\xi(t)$ and a {\it fast} partial differential equation (PDE) for the remainder function $w({\bf x},t)$ \cite{ZelMiel09,pauthier2021,rossides2023}.

In \cite{rossides2023}, the efficiencies of this numerical scheme were examined on a canonical equation for pattern formation and nonlinear optics, the one-dimensional quintic complex Ginzburg-Landau equation, which is \red{well-known} to possess spatial reflection, symmetric pulses \cite{Moores93,Akhmediev972,Mandel04,Ding2011}. Hence, the parameterization of these pulses are the shift in the $x-$direction, $r^x\in\mathbb{R}$ and $e^{\ri g}\in\mathbb{S}_1$ where $g$ is the phase shift. Under the assumption that the original one-pulse manifold is hyperbolic, then the centre-manifold reduction can be justified. In \cite{rossides2023} they focused extensively on justifying the numerical implementation of this centre-manifold approach for the interaction of two pulses. In this case the slow ODEs, in the absence of the remainder function, led to a Hamiltonian system, consisting of a phase-plane with closed cells of periodic orbits oscillating around a centre equilibrium, with heteroclinic orbits linking saddle points defining the edge of each cell. When the remainder function was included, the centre-manifold reduction approach was able to rigorously show that these heteroclinic orbits `split', and the centre equilibria become foci for the phase trajectories, which either spiral into these foci, or in some parameter scenarios, spiral into a stable limit cycle. For three pulse interactions the centre-manifold reduction approach was able to identify equilibrium solutions of the system, such as fixed points, which are useful for characterizing the existence, bifurcation, and local stability of coherent structures in infinite-dimensional dynamical systems \cite{vanderbauwhede1992,haragus2010,carr2012}.

In this paper, we extend the centre-manifold reduction approach of \cite{rossides2023} to incorporate two-space dimensions, and identify the possible dynamics of multiple pulse interactions. \red{In moving from one to two spatial dimensions, the challenges lie in the increasing size of the numerical system that needs to be solved. As we are considering symmetric pulses, the solution of the steady pulse from (\ref{eqn:disppde}) is not significantly more complicated, but the number of eigenmodes required to project the system onto the centre-manifold increases from two to three. While the size of the slow system of ODEs increases from $2N$ to $3N$, where $N$ is the number of pulses, the major challenge is that the dimension of the fast PDE for the remainder function increases from one to two dimensions. Hence we need to identify a fast and effective method to calculate this remainder function such that the computational run time of the scheme remains tractable.}

\red{The paper is laid out as follows.} The formulation of the governing fast-slow system equations is given in \S\ref{sec:formulation}, while their simplified asymptotic form in the case of well-separated pulses is given in \S\ref{sec:asymptotic}. The numerical scheme to solve the fast-slow system is given in \S\ref{sec:numerical_scheme}, while solutions to the $N=2$ problem are considered in \S\ref{sec:N2}. In \S\ref{sec:basins} we consider the steady equilibrium solutions for $N=3,4,5$ and in \S\ref{sec:coherent} we consider the interaction dynamics of two $N=3$ pulse coherent structures. Concluding remarks are given in \S\ref{sec:conclusions}.

\section{Formulation of the centre-manifold projection approach}
\label{sec:formulation}

In this investigation, we consider the interaction of multi-pulses in a two-dimensional model equation, that we choose to be the two-dimensional complex quintic Ginzburg-Landau (CQGLE) equation
\begin{equation}
u_t=\alpha \nabla^2u+\beta u+\gamma|u|^2u+\delta|u|^4u=\alpha \nabla^2 u+\beta u+ f(u).\label{eqn:gl}
\end{equation}
Here $u(x,y,t)\equiv u(r,\theta,t)$, where $r=(x^2+y^2)^{1/2}$, $\theta=\tan^{-1}(y/x)$ are plane polar coordinates with the operator
\[
\nabla^2\equiv\frac{\partial^2}{\partial x^2}+\frac{\partial^2}{\partial y^2}\equiv \frac{\partial^2}{\partial r^2}+\frac{1}{r}\frac{\partial}{\partial r}+\frac{1}{r^2}\frac{\partial^2}{\partial \theta^2},
\]
and the function
\[
f(u)=\zeta(|u|^2)u=\gamma|u|^2u+\delta|u|^4u,
\]
contains the nonlinearity of the problem, with $\alpha,\beta,\gamma$ and $\delta$ complex constants. We choose to study equation (\ref{eqn:gl}) because this is a canonical equation for pattern formation and nonlinear optics, and it is one of the simplest equations that can generate the complexity of real patterns.

The equation~\eqref{eqn:gl} has been shown to emit two-dimensional, axisymmetric, steady pulse solutions of the form $V(x,y)=V(r)$ (see for instance~\cite{Malomed2022,MALOMED2022b}), where this pulse satisfies the steady equation
\begin{equation}
\alpha\left(V_{rr}+\frac{1}{r}V_r\right)+\beta V+\gamma|V|^2V+\delta|V|^4V=0.
\label{eqn:pulse}
\end{equation}
The Ginzburg-Landau equation (\ref{eqn:gl}) can now be linearised about this steady solution (\ref{eqn:pulse}) and the resulting linear equation has a linear operator of the form
\begin{equation}
\mathbb{L}\equiv \alpha \nabla^2+\beta +f'(V),
\label{eqn:lineargl}
\end{equation}
with linearized equation
\begin{equation}
\mathbb{L}z=\alpha\nabla^2z+\beta z+\zeta(|V|^2)z+|V|^2\zeta'(|V|^2)z+V^2\zeta'(|V|^2)\overline{z}=0,
\end{equation}
where here the $\nabla^2$ operator retains the $\theta$ derivatives, which are neglected in the solution of the steady pulse, and the overbar denotes the complex conjugate.

We now consider $N$ of these axisymmetric pulses at different positions in the two-dimensional plane $(x,y)=(r_k^x,r_k^y)$ and with different phases, $g_k$, which is possible because (\ref{eqn:gl}) posses the shift symmetries $x\mapsto x-r^x$, $y\mapsto y-r^y$ and the phase symmetry $u\mapsto e^{\ri g}u$. Each of these pulses then takes the form
\begin{equation}\label{eqn:vk}
V_k(x,y,t)=V(\rb_k)=e^{\ri g_k(t)}V(x-r_k^x(t),y-r^y_k(t)),
\end{equation}
for $k=1,...,N$, where $\rb_k=(r_k^x,r_k^y,g_k)$ and these variables depend on time, $t$. In the derivation of the governing equations that follows, we will suppress the variation of these quantities with respect to $t$ for brevity. To identify a solution of (\ref{eqn:gl}) due to these $N$ pulses, we assume the solution can be written as a sum of the $N$ pulses and a remainder function, namely
\begin{equation}
u(x,y,t)=\sum_{k=1}^NV_k(x,y,t)+w(x,y,t).
\label{eqn:u}
\end{equation}

We are now interested in deriving governing differential equations for the positions and phases of each pulse, as well as the form of the remainder function, assuming that the pulses are weakly interacting, i.e. that the remainder function is small in comparison to the pulse magnitudes. To do this, we follow the global centre-manifold reduction technique laid out in \cite{rossides2023}, which will allow us to derive a system of {\it slow} ODEs for the locations and phases of the pulses, and a {\it fast} PDE for the remainder function. Substituting the solution form (\ref{eqn:u}) into (\ref{eqn:gl}) gives
\[
\sum_{k=1}^N\frac{\partial V_k}{\partial t}+\frac{\partial w}{\partial t}=\alpha\sum_{k=1}^N\left[\nabla^2 V_k+\nabla^2 w\right]+\beta\left(\sum_{k=1}^NV_k+w\right)+f\left(\sum_{k=1}^NV_k+w\right).
\]
However, from the steady solution (\ref{eqn:pulse}) we know that
\[
\alpha\nabla^2V_k=-\beta V_k-f(V_k),
\]
for each $k$, and thus
\[
\sum_{k=1}^N\frac{\partial V_k}{\partial t}+\frac{\partial w}{\partial t}=\alpha\nabla^2 w+\beta w+f\left(\sum_{k=1}^NV_k+w\right)-\sum_{k=1}^Nf(V_k).
\]
Subtracting  $f'\left(\displaystyle\sum_{k=1}^N V_k\right)w$ from each side of this equation allows us to write
\begin{equation}
\sum_{k=1}^N\frac{\partial V_k}{\partial t}+\frac{\partial w}{\partial t}-f'\left(\sum_{k=1}^N V_k\right)w 
=\alpha\nabla^2 w+\beta w+G(\Vb,w)+\Phi(\Vb),\label{eqn:goveningeqn1}
\end{equation}
where we have defined the functions
\begin{eqnarray*}
G(\Vb,w)&=&f\left(\sum_{k=1}^NV_k+w\right)-f'\left(\sum_{k=1}^N V_k\right)w-f\left(\sum_{k=1}^NV_k\right),\\
\Phi(\Vb)&=&f\left(\sum_{k=1}^NV_k\right)-\sum_{k=1}^Nf(V_k),
\end{eqnarray*}
and $\Vb=(V_1,,...,V_N)$. The function $\Phi$ is the so called {\it interaction function}, and for small $||w||$ for some suitable norm, the function $G$ is a quadratic function in $w$, which can be seen by Taylor expanding the first term, or noting that
\[
G(\Vb,0)=G_w(\Vb,0)=0.
\]

From the form of (\ref{eqn:vk}) we note that
\begin{eqnarray}
\frac{\partial V_k}{\partial t}
&=&\ri e^{\ri g_k}V(x-r_k^x,y-r_k^y)\bigdot{g}_k-e^{\ri g_k}V_x(x-r_k^x,y-r_k^y)\bigdot{r}^x_k-e^{\ri g_k}V_y(x-r_k^x,y-r_k^y)\bigdot{r}^y_k,\nonumber\\
&=&\phi_{g_k}\bigdot{g}_k+\phi_{r^x_k}\bigdot{r}^x_k+\phi_{r^y_k}\bigdot{r}^y_k,\label{eqn:vk_t}
\end{eqnarray}
where $\phi_{r^x_k}(x,y,t),\phi_{r^y_k}(x,y,t)$ and $\phi_{g_k}(x,y,t)$ are the associated zero eigenfunctions of the linear operator (\ref{eqn:lineargl}), and the overdots denote $\dfrac{\d}{\d t}$. The eigenmodes satisfy 
\begin{eqnarray*}
\phi_{g}(x,y,t)&=&\left.\frac{\partial V_k}{\partial g_k}\right|_{\rb={\bf 0}}=\ri V(x,y),\\
\phi_{r^x}(x,y,t)&=&\left.\frac{\partial V_k}{\partial r^x_k}\right|_{\rb={\bf 0}}=-V_x(x,y),\\
\phi_{r^y}(x,y,t)&=&\left.\frac{\partial V_k}{\partial r^y_k}\right|_{\rb={\bf 0}}=-V_y(x,y),
\end{eqnarray*}
where $\rb=(r^x,r^y,g)$. The main assumption we make here, and assume holds throughout this paper, is that these are the only zero eigenvalues (i.e. we do not consider Jordan blocks at zero) and the rest of the spectrum for $\mathbb{L}$ lies in the left-hand plane of the stability plane. The notations $\phi_{r_k^x}$ etc, mean we write the shifted and modulated eigenfunctions as
\begin{eqnarray}
\phi_{g_k}&=&e^{\ri g_k}\phi_{g}(x-r_k^x,y-r_k^y),\label{eqn:shiftedeig1}\\
\phi_{r_k^x}&=&e^{\ri g_k}\phi_{r_x}(x-r_k^x,y-r_k^y),\\
\phi_{r_k^y}&=&e^{\ri g_k}\phi_{r_y}(x-r_k^x,y-r_k^y).\label{eqn:shiftedeig3}
\end{eqnarray}


Inserting (\ref{eqn:vk_t}) into (\ref{eqn:goveningeqn1}) leads to the governing {\it fast} PDE for $w(x,y,t)$ as
\begin{equation}
\frac{\partial w}{\partial t}-f'\left(\sum_{k=1}^NV_k\right)w-\alpha\nabla^2 w-\beta w=-\sum_{k=1}^N\left[\phi_{r^x_k}\bigdot{r}^x_k+\phi_{r^y_k}\bigdot{r}^y_k+\phi_{g_k}\bigdot{g}_k\right]+G(\Vb,w)+\Phi(\Vb).
\label{eqn:w}
\end{equation}
To identify the governing {\it slow} ODEs for $\rb_k=(r_k^x,r_k^y,g_k)$, we project this system onto the neutral adjoint eigenfunctions of the system. The adjoint eigenfunctions are defined in the same way as the eigenfunctions, namely $\psi_{g}(x,y,t)$, $\psi_{r^x}(x,y,t)$, $\psi_{r^y}(x,y,t)$, with the spatial and phase shifted versions given the notation $\psi_{g_k}(x,y,t)$, $\psi_{r_k^x}(x,y,t)$, $\psi_{r_k^y}(x,y,t)$ as in (\ref{eqn:shiftedeig1})-(\ref{eqn:shiftedeig3}). We also introduce the $L^2(\mathbb{C})$ inner product defined such that
\[
\langle u,v \rangle=\Re\left[\int_{-\infty}^\infty\int_{-\infty}^\infty u(x,y)\overline{v(x,y)}\,\d y \d x\right],
\]
where the over bar denotes the complex conjugate. Using this inner product, the eigenfunctions and adjoint eigenfunctions satisfy the following orthogonality properties
\begin{eqnarray}
\langle\phi_{r^x},\psi_{r^x}\rangle=\langle\phi_{r^y},\psi_{r^y}\rangle=\langle\phi_{g},\psi_{g}\rangle&=&1,\label{eqn:norm1}\\
\langle\phi_{r^x},\psi_{r^y}\rangle=\langle\phi_{r^x},\psi_{g}\rangle=\langle\phi_{r^y},\psi_{r^x}\rangle=\langle\phi_{r^y},\psi_{g}\rangle=\langle\phi_{g},\psi_{r^x}\rangle=\langle\phi_{g},\psi_{r^y}\rangle&=&0.\label{eqn:norm2}
\end{eqnarray}

The projection onto the adjoint eigenfunctions of the system is constructed such that the remainder function $w(x,y,t)$ is transversal to the adjoint eigenfunctions, meaning that
\[
\langle w,\psi_{r^x_k}\rangle=\langle w,\psi_{r^y_k}\rangle=\langle w,\psi_{g_k}\rangle=0,
\]
for each $k$. The consequence of this construction is that
\begin{eqnarray*}
\frac{\d}{\d t}\langle w,\psi_{r_k^x}\rangle
&=&\Re\left[\int_{-\infty}^\infty\int_{-\infty}^\infty \left(\frac{\partial w}{\partial t}\overline{\psi}_{r_k^x}+w\frac{\partial\overline{\psi}_{r_k^x}}{\partial t}\right)\,\d y \d x\right],\\
&=&\langle w_t,\psi_{r_k^x}\rangle+\langle w,\partial_t \psi_{r_k^x}\rangle=0.
\end{eqnarray*}
However, by the definition of $\psi_{r_k^x}$,
\begin{eqnarray*}
\frac{\partial}{\partial t} \psi_{r_k^x}
&=&e^{\ri g_k(t)}\psi_{r^x}(x-r_k^x(t),y-r^y_k(t))\ri\bigdot{g}_k
-e^{\ri g_k(t)}\frac{\partial}{\partial x} \psi_{r^x}(x-r_k^x(t),y-r^y_k(t))\bigdot{r}_k^x \\
&&-e^{\ri g_k(t)}\frac{\partial}{\partial y} \psi_{r^x}(x-r_k^x(t),y-r^y_k(t))\bigdot{r}_k^y,
\end{eqnarray*} 
thus
\[
\frac{\d}{\d t}\langle w,\psi_{r_k^x}\rangle=\langle w_t,\psi_{r_k^x}\rangle-\langle \ri w,\psi_{r_k^x}\rangle\bigdot{g}_k-\langle w,\partial_x \psi_{r_k^x}\rangle\bigdot{r}_k^x-\langle w,\partial_y \psi_{r_k^x}\rangle\bigdot{r}_k^y=0.
\]
Therefore we can identify that
\begin{eqnarray*}
\langle w_t,\psi_{r_k^x}\rangle&=&\langle \ri w,\psi_{r_k^x}\rangle\bigdot{g}_k+\langle w,\partial_x \psi_{r_k^x}\rangle\bigdot{r}_k^x+\langle w,\partial_y \psi_{r_k^x}\rangle\bigdot{r}_k^y,
\end{eqnarray*}
with similar expressions valid for $\langle w_t,\psi_{r^y_k}\rangle$ and $\langle w_t,\psi_{g_k}\rangle$.

A similar simplification can be made to the $\alpha\nabla^2w$ term in (\ref{eqn:goveningeqn1}) after the projection, by noting that integration by parts shows that
\[
\langle \alpha\nabla^2 w,\psi_{r_k^x}\rangle=\langle w,\overline{\alpha}\nabla^2\psi_{r_k^x}\rangle.
\]
However, the adjoint eigenfunctions satisfy
\[
\overline{\alpha}\nabla^2\psi_{r_k^x}+\overline{\beta}\psi_{r_k^x}+\left[f'(V_k)\right]^*\psi_{r_k^x}=0,
\]
where the $*$ denotes the adjoint, thus
\begin{align*}
\langle \alpha\nabla^2 w,\psi_{r_k^x}\rangle=&-\langle w,(\overline{\beta}+\left[f'(V_k)\right]^*)\psi_{r_k^x}\rangle,\\
=&-\langle w,\overline{\beta}\psi_{r_k^x}\rangle-\langle w,\left[f'(V_k)\right]^*\psi_{r_k^x}\rangle,\\
=&-\langle \beta w,\psi_{r_k^x}\rangle-\langle f'(V_k)w,\psi_{r_k^x}\rangle,
\end{align*}
by the definition of the adjoint. Similar expressions can be found projecting onto $\psi_{r_k^y}$ and $\psi_{g_k}$.


Therefore, projecting the remainder function equation (\ref{eqn:w}) onto the adjoint eigenfunction $\psi_{r_k^x}$ for each value of $k$, gives ODEs for $r_k^x, r_k^y$ and $g_k$ of the form
%
\begin{eqnarray}
&&\sum_{j=1}^N\left[\langle\phi_{r^x_j},\psi_{r_k^x}\rangle\bigdot{r}^x_j+\langle\phi_{r^y_j},\psi_{r_k^x}\rangle\bigdot{r}^y_j+\langle\phi_{g_j},\psi_{r_k^x}\rangle\bigdot{g}_j\right]+\langle w,\partial_x \psi_{r_k^x}\rangle\bigdot{r}_k^x+\langle w,\partial_y \psi_{r_k^x}\rangle\bigdot{r}_k^y+\langle \ri w,\psi_{r_k^x}\rangle\bigdot{g}_k \nonumber\\
&&~~~~~~~~~~~~~~~~~~~~~~~~~~~~~~~~~~~~~~~~~~~~~~~~~~~~~~~~~~~~~~~~=\langle G+\Psi_k w+\Phi,\psi_{r_k^x}\rangle,
\label{eqn:ode1}
\end{eqnarray}
where
\[
\Psi_k=f'\left(\sum_{j=1}^N V_j\right)-f'(V_k),
\]
for $k=1,...,N$. Similarly, projecting onto $\psi_{r_k^y}$ and $\psi_{g_k}$ gives the ODEs
\begin{eqnarray}
&&\sum_{j=1}^N\left[\langle\phi_{r^x_j},\psi_{r_k^y}\rangle\bigdot{r}^x_j+\langle\phi_{r^y_j},\psi_{r_k^y}\rangle\bigdot{r}^y_j+\langle\phi_{g_j},\psi_{r_k^y}\rangle\bigdot{g}_j\right]+\langle w,\partial_x \psi_{r_k^y}\rangle\bigdot{r}_k^x+\langle w,\partial_y \psi_{r_k^y}\rangle\bigdot{r}_k^y+\langle \ri w,\psi_{r_k^y}\rangle\bigdot{g}_k \nonumber\\
&&~~~~~~~~~~~~~~~~~~~~~~~~~~~~~~~~~~~~~~~~~~~~~~~~~~~~~~~~~~~~~~~~=\langle G+\Psi_kw+\Phi,\psi_{r_k^y}\rangle,
\label{eqn:ode2}
\end{eqnarray}
and
\begin{eqnarray}
&&\sum_{j=1}^N\left[\langle\phi_{r^x_j},\psi_{g_k}\rangle\bigdot{r}^x_j+\langle\phi_{r^y_j},\psi_{g_k}\rangle\bigdot{r}^y_j+\langle\phi_{g_j},\psi_{g_k}\rangle\bigdot{g}_j\right]+\langle w,\partial_x \psi_{g_k}\rangle\bigdot{r}_k^x+\langle w,\partial_y \psi_{g_k}\rangle\bigdot{r}_k^y+\langle \ri w,\psi_{g_k}\rangle\bigdot{g}_k \nonumber\\
&&~~~~~~~~~~~~~~~~~~~~~~~~~~~~~~~~~~~~~~~~~~~~~~~~~~~~~~~~~~~~~~~~=\langle G+\Psi_kw+\Phi,\psi_{g_k}\rangle.
\label{eqn:ode3}
\end{eqnarray}
This gives a system of $3N$ equations, which we can write in matrix form as
\begin{equation}
{\bf C}(\rb_1,...,\rb_N,w)\bigdot{{\bf X}}={\bf F},
\label{eqn:matrix}
\end{equation}
where
\begin{eqnarray*}
{\bf X}&=&(\rb_1,...,\rb_N)^T,\\
{\bf F}&=&(\langle G+\Psi_1w+\Phi,\psi_{r_1^x}\rangle,\langle G+\Psi_1w+\Phi,\psi_{r_1^y}\rangle,\langle G+\Psi_1w+\Phi,\psi_{g_1}\rangle, \\
&&...,\langle G+\Psi_Nw+\Phi,\psi_{r_N^x}\rangle,\langle G+\Psi_Nw+\Phi,\psi_{r_N^y}\rangle,\langle G+\Psi_Nw+\Phi,\psi_{g_N}\rangle)^T.
\end{eqnarray*}
For {\it well-separated} pulses, where the norm of the remainder function $||w||$ is small in amplitude compared to the initial pulses, then the matrix ${\bf C}$ is close to the identity, and hence is invertible.

For a given remainder function $w(x,y,t)$, (\ref{eqn:matrix}) is integrated forward in time from some initial condition ${\bf X}_0=(\rb_1(0),...,\rb_N(0))^T$ to update the time dependent values in 
\begin{equation}
\bigdot{{\bf X}}={\bf C}^{-1}{\bf F}.
\label{eqn:matrix2}
\end{equation}
This gives the positions and phases, of the pulses which are then substituted into (\ref{eqn:w}) to update the remainder function $w(x,y,t)$. Equation (\ref{eqn:w}) is solved together with the boundary conditions
\begin{equation}
w(x\to\pm\infty,y,t)=w(x,y\to\pm\infty,t)=0. 
\label{eqn:w_boundary}
\end{equation}
and the initial condition 
\begin{equation}
w(x,y,0)=0. 
\label{eqn:w_initial}
\end{equation}
This initial condition is chosen because it is close to our stable solution manifold, and so we expect our numerical solution to be indistinguishably close to the solution manifold after a short transient period \cite{rossides2023}.

Under the assumption that the pulses are {\it well-separated}, we can perform an asymptotic analysis on (\ref{eqn:w}), leading to a form of the equation which is more readily solved at leading order, yet retaining sufficient accuracy to solve for the full system dynamics. 

\section{Asymptotic structure for well-separated pulses}
\label{sec:asymptotic}

The time-dependent motion of the pulses in the plane is determined by the {\it slow} ODEs (\ref{eqn:ode1})-(\ref{eqn:ode3}) together with the {\it fast} PDE (\ref{eqn:w}), with the time-dependent distance between pulses $V_j$ and $V_k$ given by 
\begin{equation}
d_{jk}(t)=\left[(r_j^x-r_k^x)^2+(r_j^y-r_k^y)^2\right]^{1/2}. 
\label{eqn:djk}
\end{equation}
The minimum separation distance for all pulses is denoted $d$, where
\begin{equation*}
d={\rm min}_{\forall j,k~j\neq k~\forall t}\{d_{jk}(t)\}.
\end{equation*}
If the pulses are well-separated, then $d$ is large, and we define the small positive parameter
\begin{equation}
\epsilon=e^{-\lambda_rd},
\label{eqn:epsilon}
\end{equation}
which is small enough for linearization. Here $\lambda=\lambda_r+\ri\lambda_i$ is a measure of the decay of the pulses at infinity and satisfies the relation
\begin{equation}
\alpha\lambda^2+\beta=0.
\label{eqn:disp}
\end{equation}
Given condition (\ref{eqn:epsilon}), it has been shown that a stable centre-manifold exists for this fast-slow system by theorem 2.1 in \cite{rossides2023}, which is proved in \cite{ZelMiel09}. In appendix \ref{appen:theorem}, we restate this theorem for completeness and proceed here with the understanding that such a stable centre-manifold exists. In the theorem itself we note that the constant $\kappa>0$ is not usually small (it is related to the spectral bound of the linear operator $\mathbb{L}$ after excluding the zero eigenvalues), and as such if we start with an initial condition sufficiently close to the manifold, i.e. the condition in (\ref{eqn:w_initial}) puts us $\epsilon$-close to the manifold, then after a short initial transient behaviour, we will be indistinguishably close to this manifold, and hence the initial condition (\ref{eqn:w_initial}) for the remainder function is justified.





With this understanding, we can remove the fast time scales by introducing the new scaled quantities
\begin{equation}
\widehat{t}=\epsilon t,~~~~\Phi=\epsilon\widehat{\Phi},~~~w=\epsilon\widehat{w},~~~\Psi_k=\epsilon\widehat{\Psi}_k,~~~G=\epsilon^2\widehat{G},
\label{eqn:scales}
\end{equation}
where the hat variables are $O(1)$ \cite{rossides2023}. Note, neither $\Phi$ nor $\Psi_k$ are $O(\epsilon)$ everywhere in the $(x,y)$-plane, but they are $O(\epsilon)$ in the regions where they interact with the adjoint eigenfunctions. Thus, the inner products can be considered of the form $\langle\Phi,\psi_{r_k^x}\rangle=O(\epsilon)$ and $\langle\Psi_k w,\psi_{g_k}\rangle=O(\epsilon^2)$ for example.

Substituting the forms (\ref{eqn:scales}) into the fast-slow system leads to an expansion of the matrix ${\bf C}$ from (\ref{eqn:matrix}) of
\[
{\bf C}(\rb_1,...,\rb_N,w)={\bf I}_{3N}+\epsilon{\bf \widehat{C}}(\rb_1,...,\rb_N,w)+O(\epsilon^2).
\]
Thus (\ref{eqn:matrix}) becomes
\[
\epsilon({\bf I}_{3N}+\epsilon\widehat{{\bf C}}){\bf X}_{\that}=\epsilon{\bf F}_1+\epsilon^2{\bf F}_2+O(\epsilon^3),
\]
where
\begin{eqnarray*}
{\bf F}_1&=&(\langle \Phihat,\psi_{r_1^x}\rangle,\langle \Phihat,\psi_{r_1^y}\rangle,\langle \Phihat,\psi_{g_1}\rangle,...,\langle \Phihat,\psi_{r_N^x}\rangle,\langle \Phihat,\psi_{r_N^y}\rangle,\langle \Phihat,\psi_{g_N}\rangle)^T,\\
{\bf F}_2&=&(\langle \Ghat+\Psihat_1\what,\psi_{r_1^x}\rangle,\langle \Ghat+\Psihat_1\what,\psi_{r_1^y}\rangle,\langle \Ghat+\Psihat_1\what,\psi_{g_1}\rangle, \\
&&...,\langle \Ghat+\Psihat_N\what,\psi_{r_N^x}\rangle,\langle \Ghat+\Psihat_N\what,\psi_{r_N^y}\rangle,\langle \Ghat+\Psihat_N\what,\psi_{g_N}\rangle)^T.
\end{eqnarray*}
Retaining only terms at $O(1)$, leads to the leading order problem
\begin{equation}
{\bf X}_{\that}={\bf F}_1.
\label{eqn:O1}
\end{equation}
We term this system the Projected ODE System (POS). Here the fast variable $\what$ does not appear, and hence the solution of (\ref{eqn:w}) is not required and can be neglected. Equation (\ref{eqn:O1}) gives a good approximation of the system dynamics for very-well-separated pulses. If the pulses are closer but still well-separated, then we can include terms at $O(\epsilon)$ to give a system which is a correction to (\ref{eqn:O1}) and thus we solve
\begin{equation}
{\bf X}_{\that}={\bf F}_1+\epsilon\left[{\bf F}_2-\widehat{{\bf C}}{\bf F}_1\right].
\label{eqn:Ow}
\end{equation}
We term this system the Projected System (PS). However, we note that the $\what$ which appears in ${\bf F}_2$ can also be expanded for small $\epsilon$, thus with the (\ref{eqn:scales}) variable scalings, (\ref{eqn:w}) becomes
\[
\epsilon^2\frac{\partial \what}{\partial\that}-\epsilon\mathbb{L}\what=-\epsilon {\bf H}\cdot{\bf X}_{\that}+\epsilon^2\widehat{G}+\epsilon\widehat{\Phi},
\]
where ${\bf H}=(\phi_{r_1^x},\phi_{r_1^y},\phi_{g_1},...,\phi_{r_N^x},\phi_{r_N^y},\phi_{g_N})^T$. Therefore retaining terms of $O(\epsilon)$ only, the leading order form of the equation for $\what$ is
\begin{equation}
\mathbb{L}\what={\bf H}\cdot{\bf X}_{\that}-\widehat{\Phi},
\label{eqn:ow}
\end{equation}
which is a stationary PDE, and hence is simpler to solve than the full time-dependent equation (\ref{eqn:goveningeqn1}).

\subsection{Approximation of POS for very-well-separated pulses}
\label{sec:leadingorder}

At leading order the governing equations do not depend on the remainder function $w=\epsilon\what$ and so the forms of the inner products in ${\bf F}_1$ can be constructed explicitly, making their numerical evaluation fast, in the very-well-separated pulse limit. Here we construct the form of these inner products for the interaction of pulses $V_j$ and $V_k$.

The inner products in the vector ${\bf F}_1$ can be simplified by using Lemma 10.2 of \cite{ZelMiel09} whereby we can write
\begin{eqnarray}
\langle \Phi, \psi_{r^x_k} \rangle &=& -\sum_{\substack{j=1\\j\neq k}}^N \langle \mathbb{A} V_j, \psi_{r^x_k} \rangle = -\sum_{\substack{j=1\\j\neq k}}^N \langle V_j, \overline{\mathbb{A}}\psi_{r^x_k} \rangle,\label{eqn:A}\\
\langle \Phi, \psi_{r^y_k} \rangle &=& -\sum_{\substack{j=1\\j\neq k}}^N \langle \mathbb{A} V_j, \psi_{r^y_k} \rangle = -\sum_{\substack{j=1\\j\neq k}}^N \langle V_j, \overline{\mathbb{A}}\psi_{r^y_k} \rangle,\label{eqn:B}\\
\langle \Phi, \psi_{g_k} \rangle &=& -\sum_{\substack{j=1\\j\neq k}}^N \langle \mathbb{A} V_j, \psi_{g_k} \rangle = -\sum_{\substack{j=1\\j\neq k}}^N \langle V_j, \overline{\mathbb{A}}\psi_{g_k} \rangle,\label{eqn:C}
\end{eqnarray}
where $\mathbb{A}$ is the elliptic differential operator of (\ref{eqn:gl}), namely
\[
\mathbb{A}\equiv\alpha\nabla^2 +\beta.
\]

Each of the inner products can be approximated by performing the two-dimensional analogue of the one-dimensional calculation in \cite{ZelMiel09}. The calculation assumes the large $r=(x^2+y^2)^{1/2}$ asymptotic form of the pulse and the adjoint eigenmodes, and evaluates the inner product in this overlap region. It will also be useful to have the integral identity for the linear operator $\mathbb{A}$, where for functions $\xi_1(r,\theta)$ and $\xi_2(r,\theta)$
\[
\int_0^{2\pi}\int_0^\infty \xi_1(\mathbb{A}\xi_2) r\, \d r\d\theta=\alpha\int_0^{2\pi}\left[r\xi_1\xi_{2r}-r\xi_{1r}\xi_2\right]_0^\infty \d\theta+\int_0^{2\pi}\int_0^\infty (\mathbb{A}\xi_1)\xi_2 r\, \d r\d\theta,
\]
after integrating by parts twice. Therefore, for pulses at $(x,y)=(r_j^x,r_j^y)$ and $(r_k^x,r_k^y)$ we have
\small
\begin{eqnarray*}
-\langle V_j, \overline{\mathbb{A}}\psi_{g_k} \rangle&=&\int_{-\infty}^\infty\int_{-\infty}^\infty V(x-r_j^x,y-r_j^y)\mathbb{A}\overline{\psi}_g(x-r_k^x,y-r_k^y)\,\d x\d y,\\
&=&\int_{0}^{2\pi}\int_{0}^\infty V(r_k^x-r_j^x+R\cos\theta,r_k^y-r_j^y+R\sin\theta)\mathbb{A}\overline{\psi}_g(R\cos\theta,R\sin\theta) R\,\d R\d\theta,\\
&=&\int_{0}^{2\pi}\int_{0}^\infty V(r_k^x-r_j^x+R\cos\theta,r_k^y-r_j^y+R\sin\theta)\mathbb{A}\Big[\overline{\psi}_g(R\cos\theta,R\sin\theta) 
-\overline{\psi}_g^a(R\cos\theta,R\sin\theta)\Big] R\,\d R\d\theta,\\
&=&\int_0^{2\pi}\alpha\lim_{R\to0}\left[RV(r_j^x-r_k^x+R\cos\theta,r_j^y-r_k^y+R\sin\theta)\Big[\overline{\psi}_g(R\cos\theta,R\sin\theta) \right. 
-\overline{\psi}_g^a(R\cos\theta,R\sin\theta)\Big]_R\\&& \left.-RV_R(r_j-r_k+R\cos\theta,R\sin\theta)\Big[\overline{\psi}_g(R\cos\theta,R\sin\theta)-\overline{\psi}_g^a(R\cos\theta,R\sin\theta)\Big]\right] \d\theta\\
&&+\int_{0}^{2\pi}\int_{0}^\infty \mathbb{A}V(r_j-r_k+R\cos\theta,R\sin\theta)\Big[\overline{\psi}_g(R\cos\theta,R\sin\theta)-\overline{\psi}_g^a(R\cos\theta,R\sin\theta)\Big] R\,\d R\d\theta,
\end{eqnarray*}
\normalsize
where we have used $(x,y)=(r_k^x,r_k^y)+R(\cos\theta,\sin\theta)$ in the second line, and the `$a$' superscript denotes the large-$r$ asymptotic form of the adjoint eigenmode. The final integral is now small, due to the fact that the asymptotic result has been subtracted from the full eigenmode, and we have used the fact that as $R\to\infty$ the pulse and eigenmode are exponentially small. Thus
\begin{eqnarray}
-\langle V_j, \overline{\mathbb{A}}\psi_{g_k} \rangle&=&\int_0^{2\pi}\alpha\lim_{R\to0}\left[RV(r_k^x-r_j^x+R\cos\theta,r_k^y-r_j^y+R\sin\theta)
\Big[\overline{\psi}_g(R\cos\theta,R\sin\theta)-\overline{\psi}_g^a(R\cos\theta,R\sin\theta)\Big]_R \right.\nonumber\\
&& \left.-RV_R(r_k-r_j+R\cos\theta,R\sin\theta) 
\Big[\overline{\psi}_g(R\cos\theta,R\sin\theta)-\overline{\psi}_g^a(R\cos\theta,R\sin\theta)\Big]\right] \d\theta. \label{eqn:vjinner}
\end{eqnarray}

The asymptotic forms of the pulses and the adjoint eigenmodes can be determined from the governing linear equations and are 
\begin{eqnarray*}
V^a(r,\theta)&=&pH_0(\ri\lambda r),\\
\psi_g^a(r,\theta)&=&s\overline{H}_0(\ri\lambda r),\\
\psi_{r^x}^a(r,\theta)&=&q\overline{H}_0(\ri\lambda r)\cos\theta,\\
\psi_{r^y}^a(r,\theta)&=&q\overline{H}_0(\ri\lambda r)\sin\theta,
\end{eqnarray*}
where $p,~q,~s$ are complex constants, $\lambda$ is found in (\ref{eqn:disp}) and $H_n(X)=J_n(X)+\ri Y_n(X)$ is the Hankel function of the first kind of order $n$, with $J_n(X)$ and $Y_n(X)$ Bessel functions of the first kind.

As $R\to0$ in (\ref{eqn:vjinner}) then
\begin{eqnarray*}
V(r_k^x-r_j^x+R\cos\theta,r_k^y-r_j^y+R\sin\theta)&=&pH_0(\ri\lambda d_{jk})-p\lambda H_1(\ri\lambda d_{jk})R\cos\theta
-p\lambda H_1(\ri\lambda d_{jk})R\sin\theta+O(R^2),\\
\overline{\psi}_{g}^a(R\cos\theta,R\sin\theta)&=&\overline{s}\frac{2\ri}{\pi}\log(\ri\lambda R)+O(1),\\
\overline{\psi}_{g}(R\cos\theta,R\sin\theta)&=&C_1+O(R),
\end{eqnarray*}
where $C_1$ is an unimportant constant and $d_{jk}$ from (\ref{eqn:djk}) is the separation distance of the two pulses. Therefore, as $R\to0$ the inner product using (\ref{eqn:C}) becomes
\[
\langle \Phi, \psi_{g_k} \rangle=-\sum_{\substack{j=1\\j\neq k}}^N 4\ri\alpha p\overline{s}H_0(\ri\lambda d_{jk})e^{-\ri(g_k-g_j)}.
\]

For the inner products involving $\psi_{r^x}$ and $\psi_{r^y}$ there is a problem to overcome, as
\[
\langle \Phi, \psi_{r_k^x} \rangle=\langle \Phi, \psi_{r_k^y} \rangle=O(R\log(R)),
\]
as $R\to0$ and hence would give a zero inner product, whereas numerical evaluations give these as a non-zero inner product. This issue is down to the singular form of the asymptotic eigenfunctions for the 2D problem, but this issue is not evident in the 1D problem in \cite{ZelMiel09} where the asymptotic eigenmodes are just exponential functions. Examining the above limit as $R\to0$, and considering the one-dimensional result in \cite{ZelMiel09}, it can be seen that $\langle \Phi,\psi_{r_k^x}\rangle$ and $\langle \Phi,\psi_{r_k^y}\rangle$ have a form
\[
\langle \Phi, \psi_{r^x_k} \rangle=-\sum_{\substack{j=1\\j\neq k}}^N B\alpha p\overline{q}\pi\lambda \frac{(r_k^x-r_j^x)}{d_{jk}}H_1(\ri\lambda d_{jk})e^{-\ri(g_k-g_j)},
\]
and
\[
\langle \Phi, \psi_{r^y_k} \rangle=-\sum_{\substack{j=1\\j\neq k}}^N B\alpha p\overline{q}\pi\lambda \frac{(r_k^y-r_j^y)}{d_{jk}} H_1(\ri\lambda d_{jk})e^{\ri(g_k-g_j)},
\]
where $B$ is a complex constant. This asymptotic result will be confirmed numerically in \S\ref{sec:numerical_scheme}, where the value of $B$ will be estimated.

In the case of two pulses, which we consider as both being located on the $x$-axis without loss of generality, we can write the four ODE system as
\[
\bigdot{r}_1^x=\langle \Phi, \psi_{r^x_1} \rangle,~~~
\bigdot{r}_2^x=\langle \Phi, \psi_{r^x_2} \rangle,~~~
\bigdot{g}_1=\langle \Phi, \psi_{g_1} \rangle,~~~
\bigdot{g}_2=\langle \Phi, \psi_{g_2} \rangle.
\]
In this case we are able to simplify the system to a pair of ODEs for the difference functions $\rbar=r_1^x-r_2^x$, $\gbar=g_1-g_2$, which by using the following asymptotic forms of the Hankel functions
\begin{eqnarray}
H_0(\ri\lambda\rbar)&=&-\sqrt{\frac{2\ri}{\pi\lambda|\rbar|}}e^{-\lambda_r|\rbar|}e^{\ri(-\lambda_i|\rbar|+\pi/4)}+O\left(e^{-\lambda_r|\rbar|}|\rbar|^{-3/2}\right),\label{eqn:H0}\\
H_1(\ri\lambda\rbar)&=&\ri\sqrt{\frac{2\ri}{\pi\lambda|\rbar|}}e^{-\lambda_r|\rbar|}e^{\ri(-\lambda_i|\rbar|+\pi/4)}+O\left(e^{-\lambda_r|\rbar|}|\rbar|^{-3/2}\right), \label{eqn:H1}
\end{eqnarray}
as $|\rbar|\to\infty$, can be written as
\begin{eqnarray}
\bigdot{\rbar}&=&2J_1|\rbar|^{-1/2}e^{-\lambda_r|\rbar|}\cos(-\lambda_i|\rbar|+\kappa_1+\pi/4)\cos(\gbar),\label{eqn:rbar2}\\
\bigdot{\gbar}&=&-2J_2|\rbar|^{-1/2}e^{-\lambda_r|\rbar|}\sin(-\lambda_i|\rbar|+\kappa_2+\pi/4)\sin(\gbar),\label{eqn:gbar2}
\end{eqnarray}
where
\begin{equation}
J_1e^{\ri\kappa_1}=\ri B\alpha p\overline{q}\sqrt{2\lambda\pi}~~~~J_2e^{\ri\kappa_2}=-4\alpha p\overline{s}\sqrt{\frac{2}{\lambda\pi}}.
\label{eqn:consts}
\end{equation}

\begin{figure}[!htb]
\begin{center}
(a)\includegraphics[width=0.45\textwidth]{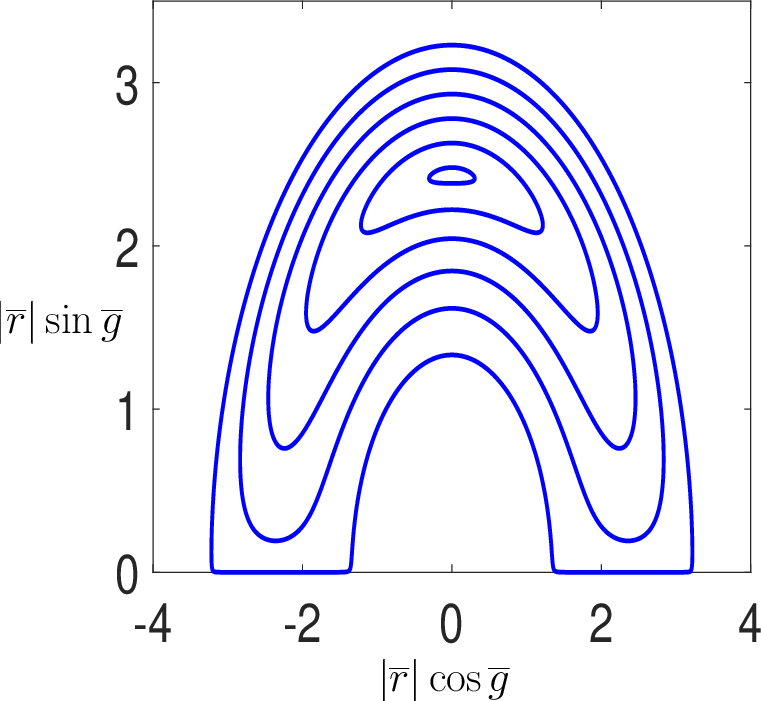}
(b)\includegraphics[width=0.45\textwidth]{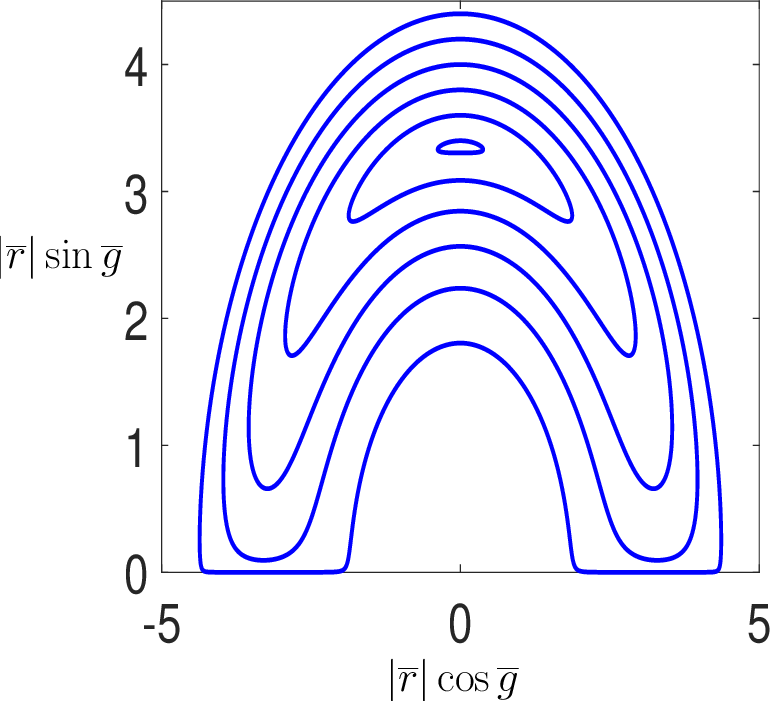}
\end{center}
\caption{Plot of the cell 1 phase space for (\ref{eqn:rbar2}) and (\ref{eqn:gbar2}) in the $(\overline{r}\cos \overline{g},\overline{r}\sin \overline{g})$-plane for the parameters in (a) (\ref{eqn:parameters1}) and (b) (\ref{eqn:parameters2}).}
\label{fig:phase_space}
\end{figure}

The form of equations (\ref{eqn:rbar2})-(\ref{eqn:gbar2}) have the same structure as those found for the one-dimensional problem in \cite{ZelMiel09,rossides2023}, and depending on the sign of $J_1J_2\lambda_i\cos(\kappa_2-\kappa_1)$, this system can exhibit a Hamiltonian structure with Hamiltonian
\[
\mathscr{H}=\sin(\gbar)\exp\left[-\frac{J_2}{J_1}\rbar\sin(\kappa_2-\kappa_1)\right]\left|\sin\left(-\lambda_i\rbar+\kappa_1+\frac{\pi}{4}\right)\right|^{\frac{J_2}{J_1\lambda_i}\cos(\kappa_2-\kappa_1)},
\]
which has the same structure as for the one-dimensional problem in \cite{Turaev07}. Our interest is in the case where $J_1J_2\lambda_i\cos(\kappa_2-\kappa_1)>0$ which leads to Hamiltonian (reversible) dynamics. The phase plane for this two pulse problem is given in Figure \ref{fig:phase_space}, and the contours are closed with equilibria consisting of centres and saddle points. The centres correspond to equilibrium points of (\ref{eqn:rbar2}) and (\ref{eqn:gbar2}) with $\cos\gbar=0$ ($\sin\gbar=\pm1$) and hence $\rbar$ satisfies
\[
\sin\left(-\lambda_i|\rbar|+\kappa_2+\frac{\pi}{4}\right)=0,
\]
thus
\[
|\rbar|=\frac{\kappa_2+\frac{\pi}{4}+n\pi}{\lambda_i},~~~~n=1,2,3,4,....
\]
Similarly, the saddles are equilibrium points with $\sin\gbar=0$ ($\cos\gbar=\pm1$) and hence
\[
|\rbar|=\frac{\kappa_1+\frac{\pi}{4}+(2n-1)\frac{\pi}{2}}{\lambda_i},~~~~n=1,2,3,4,....
\]
The saddle points are connected via heteroclinic orbits which split the phase plane into different cells, each containing periodic orbits. The cell shown in each panel of Figure \ref{fig:phase_space} are the first such cell or the cell 1 dynamics.

In one-space dimension, \cite{rossides2023} showed that this Hamiltonian structure breaks down when the $O(\epsilon)$ terms (i.e. incorporating $\widehat{w}$) are included in the analysis. It was shown that the first cell contained a stable focus or stable limit cycle, depending upon the values of $\alpha,~\beta,~\gamma,~\delta$, and all initial conditions converged to these stable results, including those initial conditions inside cells 2, 3 etc. We expect a similar structure of solution here for the two-dimensional CQGLE, which we examine in the next section.

\section{Results}

\subsection{Numerical Scheme}
\label{sec:numerical_scheme}

The steady localised pulse solutions we are seeking in this work are axisymmetric, and thus we can numerically solve for the pulse by solving the axisymmetric form of the steady CQGLE in (\ref{eqn:pulse}). We set up this problem as a boundary value problem for $V(r)$ on the finite truncated domain $r\in[0,L]$. To solve this equation, we split the solution into real and imaginary parts, $V=V_r+\ri V_i$, and solve directly the resulting coupled set of equations. The axisymmetry of the problem  means we solve these equations subject to the origin conditions
\begin{equation}
\left.\frac{\d V_r}{\d r}\right|_{r=0}=\left.\frac{\d V_i}{\d r}\right|_{r=0}=0,
\label{eqn:pulse_origin}
\end{equation}
while at $r=L$ we require the correct exponential decay at infinity, hence we impose
\begin{equation}
\frac{\d V_r}{\d r}+\lambda_r V_r-\lambda_i V_i=\frac{\d V_i}{\d r}+\lambda_i V_r+\lambda_r V_i=0,
\label{eqn:pulse_L}
\end{equation}
where $\lambda=\lambda_r+\ri\lambda_i$ satisfies (\ref{eqn:disp}).

The problem has 8 real parameters from the real and imaginary parts of $\alpha,~\beta,~\gamma,~\delta$. We fix 7 of these and find the eighth, $\beta_i={\rm Im}(\beta)$, by imposing (\ref{eqn:disp}). Equation (\ref{eqn:pulse}) along with (\ref{eqn:pulse_origin}), (\ref{eqn:pulse_L}) and (\ref{eqn:disp}) are solved using the \lstinline[style=Matlab-editor]{chebop} routine of chebfun for nonlinear boundary value problems, with a tolerance of machine precision \cite{Driscoll2014}.

To calculate the eigenfunctions and adjoint eigenfunctions of the zero eigenvalues, we consider a solution of the form $\Vhat(r)e^{\ri m\theta}$ to the linear and adjoint problems $\mathbb{L}(\Vhat e^{\ri m\theta})=\mathbb{L}^*(\Vhat e^{\ri m\theta})=0$. Here there are three zero eigenfunctions (as opposed to two for the one-dimensional problem), one with $m=0$ and two with $m=1$. With $m$ fixed at $0,1$, the linear and adjoint problems are again solved by considering real and imaginary parts $\Vhat=\Vhat_r+\ri\Vhat_i$, and again using the same chebfun structure as before, but now the origin and boundary conditions
\[
\left.\frac{\d\Vhat_r}{\d r}\right|_{r=0}=\left.\frac{\d\Vhat_i}{\d r}\right|_{r=0}=\Vhat_r(L)=\Vhat_i(L)=0,
\]
for the $m=0$ eigenfunction and 
\[
\Vhat_r(0)=\Vhat(0)=\Vhat_r(L)=\Vhat_i(L)=0,
\]
for the $m=1$ eigenfunctions. For $m=1$, the two eigenfunctions come from the real and imaginary parts of $\Vhat e^{\ri\theta}$. The resulting zero and adjoint eigenfunctions then need to be normalized for the projected system such that (\ref{eqn:norm1}) and (\ref{eqn:norm2}) hold.

The projection scheme requires translations and phase variations of the pulses and eigenfunctions, such as in (\ref{eqn:vk}). This is achieved numerically by firstly projecting all these functions back to the Cartesian domain $(x,y)\in[-L,L]\times[-L,L]$, and then, because the functions are constructed using Chebyshev polynomials in \verb1chebfun1, the translations are calculated via spectral interpolation within \verb1chebfun1.

The time-evolution for $(\rb_1,...,\rb_N)^T$ can now be determined by solving (\ref{eqn:ode1})-(\ref{eqn:ode3}) with the given initial conditions $(\rb_1(0),...,\rb_N(0))^T$ and $\what(x,y,0)=0$, and integrating forward in time using \lstinline[style=Matlab-editor]{MATLAB}'s \lstinline[style=Matlab-editor]{ode15s} or \lstinline[style=Matlab-editor]{ode45} solvers. Both \lstinline[style=Matlab-editor]{ode15s} and \lstinline[style=Matlab-editor]{ode45} are variable size time-stepping routines, and their accuracy is determined by the function and time-step errors which are both set equal to ${\rm Err}=10^{-8}$. To update $\widehat{w}$ at each time-step we solve (\ref{eqn:ow}) by discretizing the domain $(x,y)\in[-L,L]\times[-L,L]$ on an $M\times M$ grid, and making use of \lstinline[style=Matlab-editor]{MATLAB}'s sparse matrix storage ability. The resulting large system of nonlinear equations are then solved via \lstinline[style=Matlab-editor]{gmres} with the boundary conditions (\ref{eqn:w_boundary}). The inner products in (\ref{eqn:ode1})-(\ref{eqn:ode3}) are calculated using Boole's rule in two-dimensions.

\begin{figure}[!htb]
\begin{center}
(a)\includegraphics[width=0.40\textwidth]{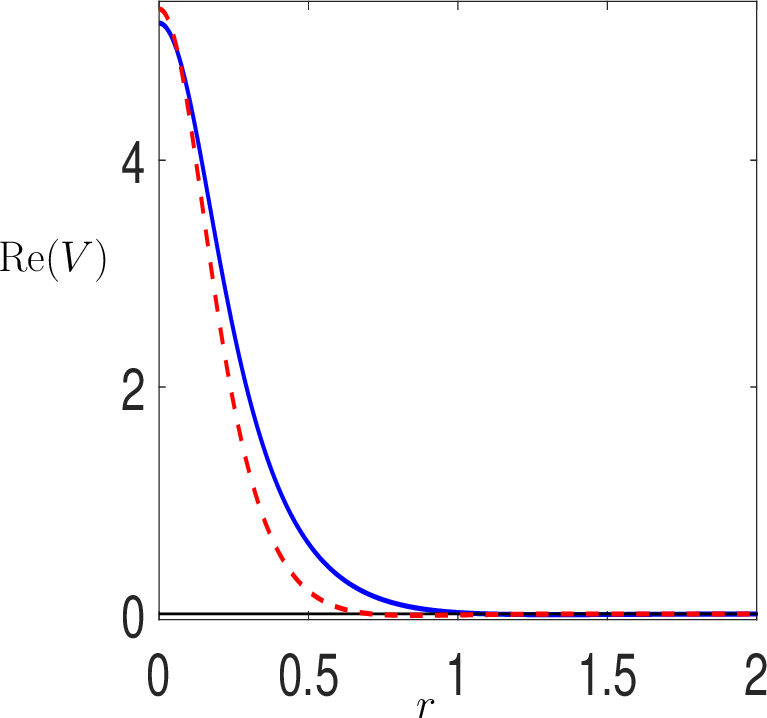}
(b)\includegraphics[width=0.40\textwidth]{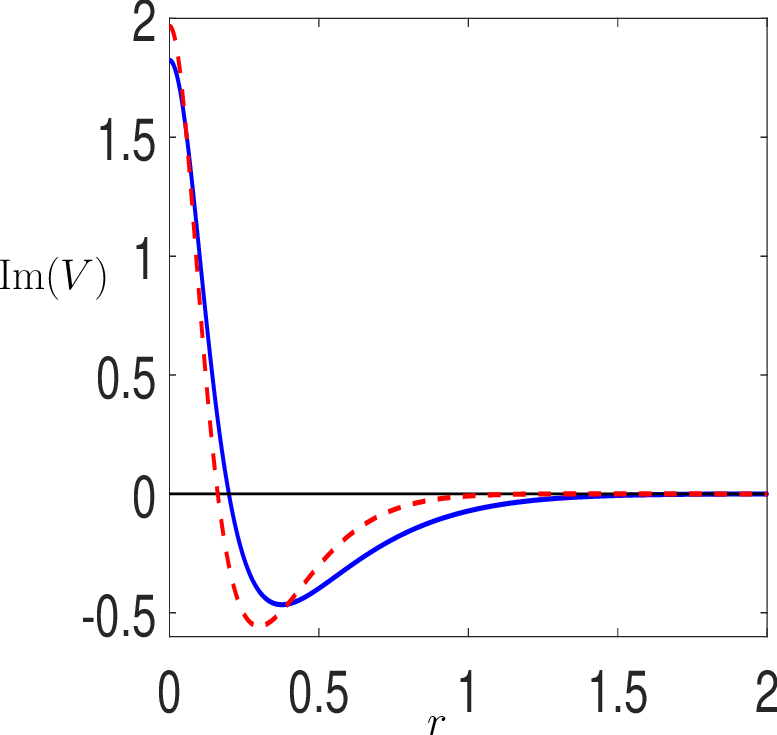}
(c)\includegraphics[width=0.40\textwidth]{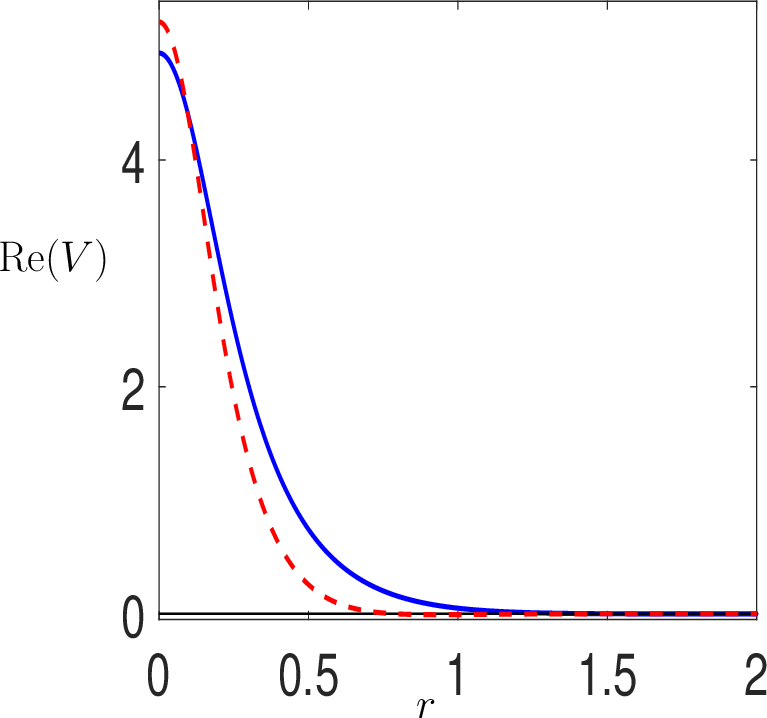}
(d)\includegraphics[width=0.40\textwidth]{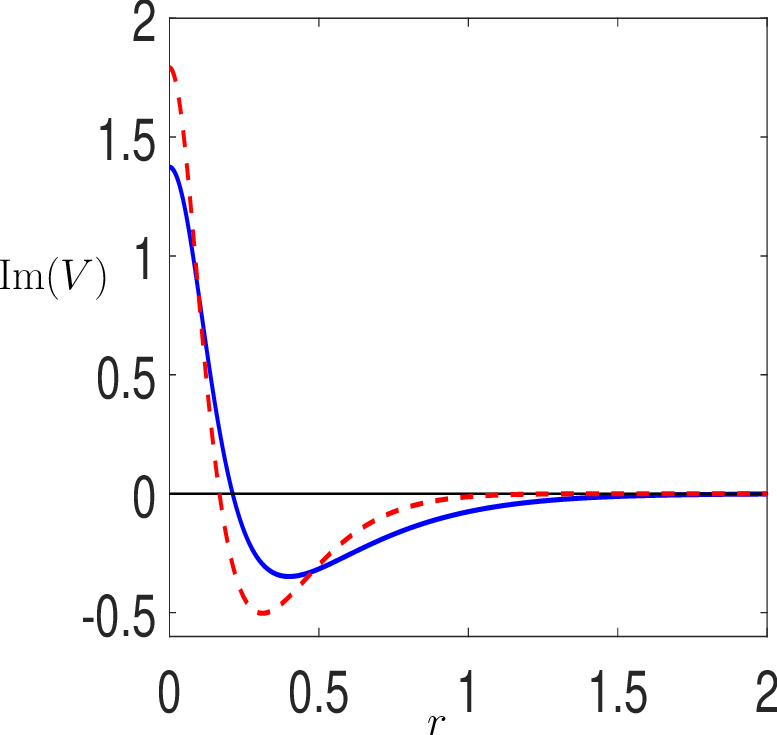}
\end{center}
\caption{Plot of (a,c) ${\rm Re}(V)$ and (b,d) ${\rm Im}(V)$ panels (b,d) for the parameter values  (\ref{eqn:parameters1}) in panels (a,b) and (\ref{eqn:parameters2}) in panels (c,d). In each plot, the solid blue line represents the two-dimensional result while the red dashed line give the corresponding one-dimensional result.}
\label{fig:pulse}
\end{figure}

For the results presented in this paper, we consider two specific axisymmetric pulse solutions with parameters
\begin{equation}
\alpha=0.5+0.5\ri, \quad \beta=-0.05-13.2\ri, \quad \gamma=1.8+\ri, \quad \delta=-0.05+0.05\ri,
\label{eqn:parameters1}
\end{equation} 
and
\begin{equation}
\alpha=0.5+0.5\ri, \quad \beta=-2.0-10.8\ri, \quad \gamma=1.8+\ri, \quad \delta=-0.05+0.05\ri,
\label{eqn:parameters2}
\end{equation} 
where the value of $\beta_i$ in each case was found as part of the initial pulse calculation. The real and imaginary parts of these pulses are plotted in Figure \ref{fig:pulse}, together with their one-dimensional equivalent, where in the one-dimensional case $\beta=-0.05-37.9\ri$ and $\beta=-2.0-34.7\ri$ respectively.

\begin{figure}[!htb]
\begin{center}
(a)\includegraphics[width=0.45\textwidth]{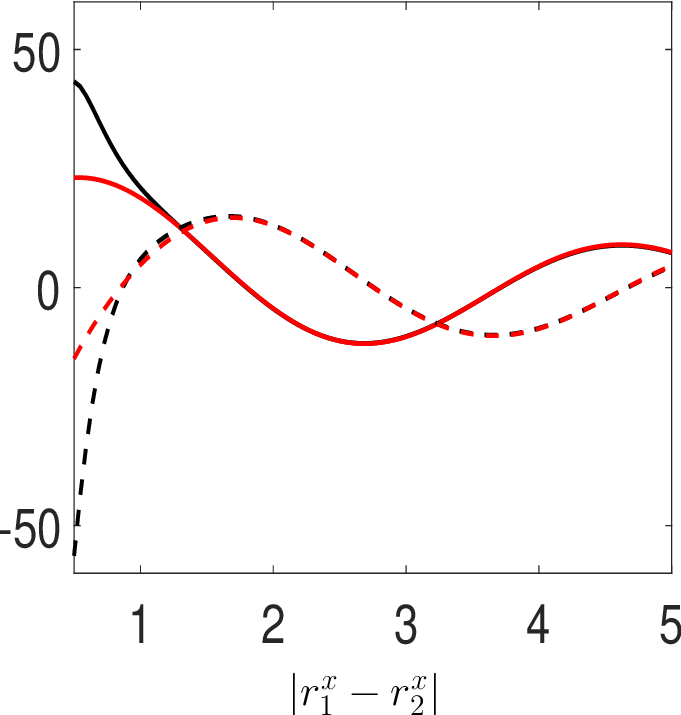}
(b)\includegraphics[width=0.45\textwidth]{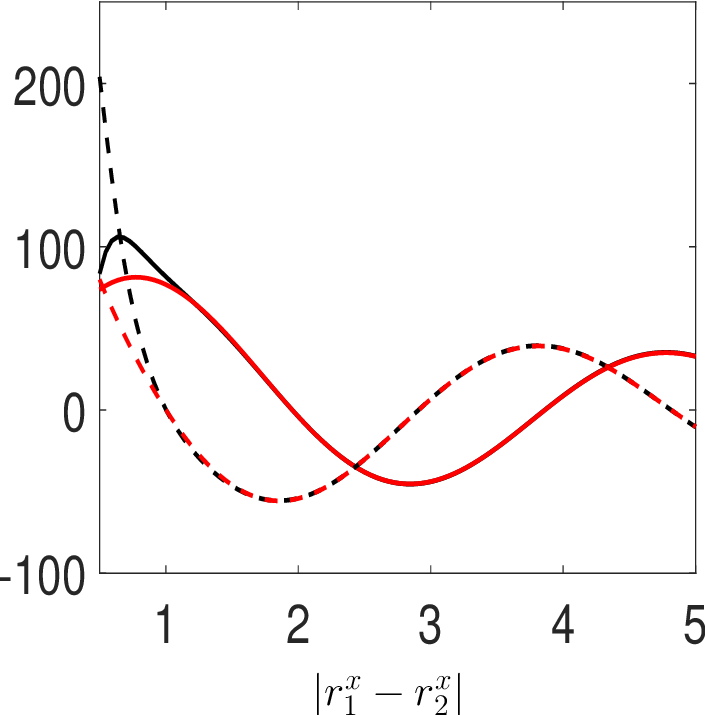}
\end{center}
\caption{(a) $\langle\Phi,\psi_{r_1^x}\rangle$ (solid lines) and $\langle\Phi,\psi_{r_2^x}\rangle$ (dashed lines) and (b) $\langle\Phi,\psi_{g_1}\rangle$ (solid lines) and $\langle\Phi,\psi_{g_2}\rangle$ (dashed lines) as a function of $|r_1^x-r_2^x|$. The black result signify the full numerical results and the red results signify the asymptotic result. In each case the inner product is multiplied by $e^{\lambda_i d_{12}}$ to remove the exponential decay of the inner product with pulse separation distance.}
\label{fig:inner_prod}
\end{figure}
In \S\ref{sec:asymptotic} we identified the leading order asymptotic forms of the pulse interaction terms, and here we are able to verify these forms numerically and provide a numerical calculation to determine the value of the constants in (\ref{eqn:consts}) for the pulse parameters given in (\ref{eqn:parameters1}). In order to do this, we consider $N=2$ pulses such that when using the large-$\rbar$ form of the Hankel functions (\ref{eqn:H0})-(\ref{eqn:H1}) then
\begin{eqnarray}
\langle \Phi,\psi_{r_1^x}\rangle&=&J_1e^{\ri\kappa_1}\frac{(r_1^x-r_2^x)}{d_{12}}e^{\ri(g_1-g_2)}e^{-\lambda_rd_{12}}e^{\ri(-\lambda_id_{12}+\pi/4)},\label{eqn:ip1}\\
\langle \Phi,\psi_{r_2^x}\rangle&=&J_1e^{\ri\kappa_1}\frac{(r_2^x-r_1^x)}{d_{12}}e^{\ri(g_2-g_1)}e^{-\lambda_rd_{12}}e^{\ri(-\lambda_id_{12}+\pi/4)},\\
\langle \Phi,\psi_{r_1^y}\rangle&=&J_1e^{\ri\kappa_1}\frac{(r_1^y-r_2^y)}{d_{12}}e^{\ri(g_1-g_2)}e^{-\lambda_rd_{12}}e^{\ri(-\lambda_id_{12}+\pi/4)},\\
\langle \Phi,\psi_{r_2^y}\rangle&=&J_1e^{\ri\kappa_1}\frac{(r_2^y-r_1^y)}{d_{12}}e^{\ri(g_2-g_1)}e^{-\lambda_rd_{12}}e^{\ri(-\lambda_id_{12}+\pi/4)},\\
\langle \Phi,\psi_{g_1}\rangle&=&J_2e^{\ri\kappa_2}e^{\ri(g_1-g_2)}e^{-\lambda_rd_{12}}e^{\ri(-\lambda_id_{12}+\pi/4)},\\
\langle \Phi,\psi_{g_2}\rangle&=&J_2e^{\ri\kappa_2}e^{\ri(g_2-g_1)}e^{-\lambda_rd_{12}}e^{\ri(-\lambda_id_{12}+\pi/4)},\label{eqn:ip6}
\end{eqnarray}
where we have omitted the $O\left(e^{-\lambda_rd_{12}}d_{12}^{-3/2}\right)$ correction terms for brevity. The unknown real coefficients $J_1,~J_2,~\kappa_1,~\kappa_2$ can then be calculated for the parameters in (\ref{eqn:parameters1}) by calculating the inner products for various values of $r_1^x-r_2^x$ and $g_1-g_2$ and approximating the best fit surface using \lstinline[style=Matlab-editor]{MATLAB}'s \lstinline[style=Matlab-editor]{fit} subroutine. Note, we simplify the calculation by setting $r_1^y-r_2^y=0$ throughout the fitting procedure, which we are able to do by choosing the pulses to lie on the $x$-axis. The \lstinline[style=Matlab-editor]{fit} subroutine gives a $95\%$ confidence interval for the complex coefficients $J_1e^{\ri\kappa_1},~J_2e^{\ri\kappa_2}$ for a numerically specified function. To make sure we are in the very-weakly-interacting regimes, we choose to evaluate the inner products for values
\[
r_1^x-r_2^x\in[2.2,5.7],~~~~~~g_1-g_2\in[-\pi,\pi],
\]
to perform the fit. Taking the central value in the confidence interval, we find
\begin{eqnarray*}
J_1&=&19.43732~~~~~~~~~~\kappa_1=-0.17263,\\
J_2&=&76.88165~~~~~~~~~~\kappa_2=~0.07856,
\end{eqnarray*}
to 5 d.p.. In Figure \ref{fig:inner_prod}, we plot the full numeric and asymptotic forms of the inner products (\ref{eqn:ip1})-(\ref{eqn:ip6}) (multiplied by $e^{\lambda_rd_{12}})$ as a function of $|d_{12}|=|r_1^x-r_2^x|$, where the pulses are located at
\[
(r_1^x,r_1^y,g_1)=(-|d_{12}|/2,0,0),~~~~~~(r_2^x,r_2^y,g_2)=(|d_{12}|/2,0,\pi/4).
\]
Here we see that the approximation and the full numerical result agree well for $|r_1^x-r_2^x|\gtrsim1.7$. For comparison, the approximation in (\ref{eqn:consts}) gives
\[
J_2=76.71014,~~~~~~~~~~\kappa_2=0.07910,
\]
which is within $1\%$ of the value found above by the numerical fitting procedure. Also note, comparing the approximation found in (\ref{eqn:consts}) for $J_1e^{\kappa_1}$ with the numerical results above gives
\[
B\approx 0.277-0.117\ri.
\]

\subsection{Case $N=2$ pulse interactions}
\label{sec:N2}

The case of $N=2$ pulses in one-space dimension was extensively studied in \cite{rossides2023}, and as the two pulse problem in the plane effectively reduces to a quasi-one-dimensional problem, we expect similarities with our results presented here. However, there are subtle differences between the one and two-space dimension problems, namely the shape and far-field decay rates of the pulses, see Figure \ref{fig:pulse}, and so we examine the two-dimensional analogue of this problem in order to validate our numerical scheme.

As noted in \S\ref{sec:leadingorder}, the leading order two pulse problem is Hamiltonian and hence the phase-plane portraits in the $(|\rbar|\sin\gbar,|\rbar|\cos\gbar)$-plane (Figure \ref{fig:phase_space}) consist of closed periodic orbits, where $\rbar=r_1^x-r_2^x$ and $\gbar=g_1-g_2$. Note here we have chosen to fix $r_1^y=r_2^y=0$, but the phase plane is identical for any $r_1^y$ and $r_2^y$ values, except now the pulse interaction occurs along a line in the $(x,y)$-plane passing through the centre of the two pulses, rather than along the $x$-axis, as we have constructed here. Figure \ref{fig:phase_space} shows the first cell of these periodic orbits with additional, weaker interacting cells, existing for larger $|\rbar|$ values.

\begin{figure}[!htb]
\begin{center}
(a)\includegraphics[width=0.4\textwidth]{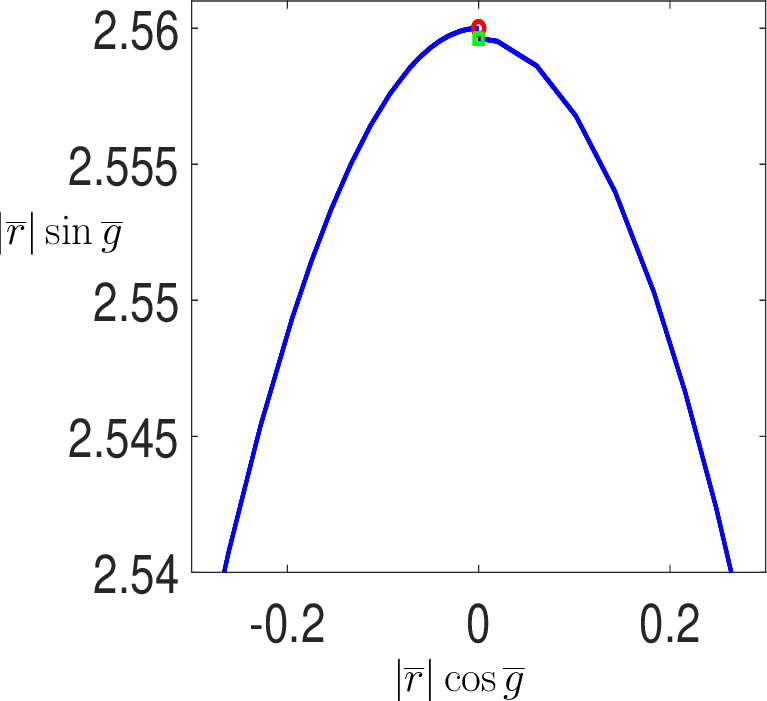}
(b)\includegraphics[width=0.4\textwidth]{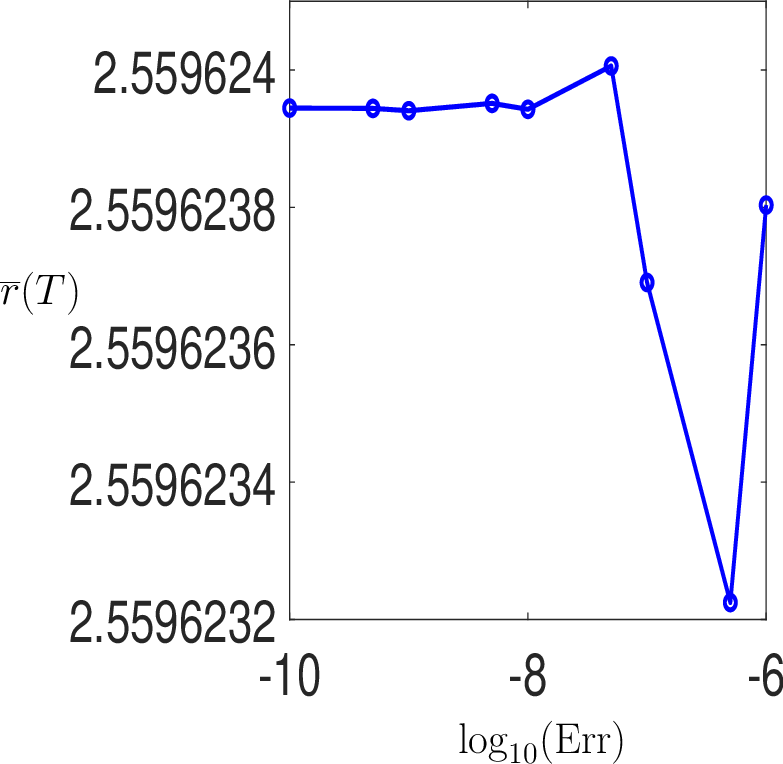}
(c)\includegraphics[width=0.4\textwidth]{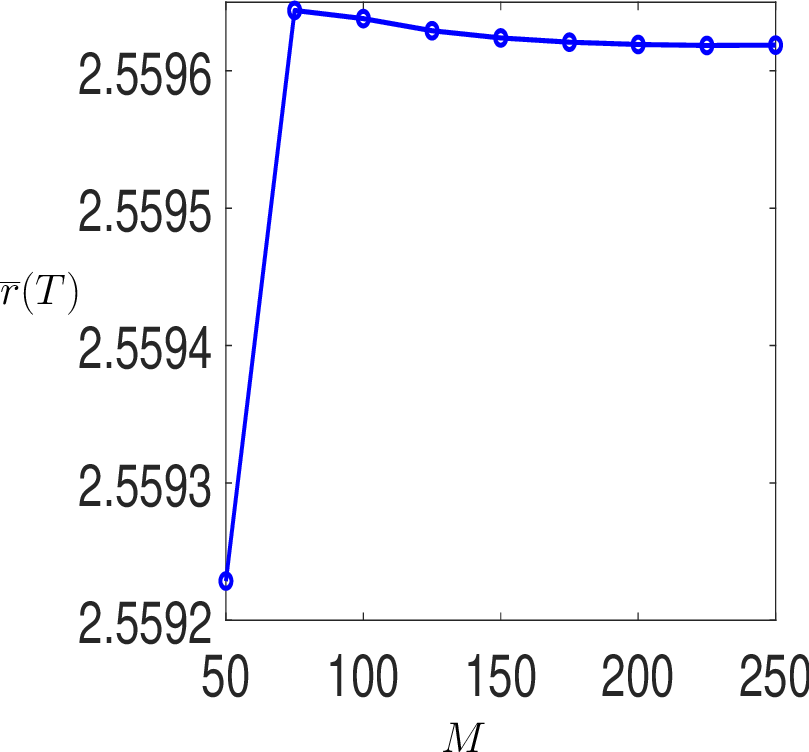}
\end{center}
\caption{(a) Blow up of the $\rbar(0)=2.56$ phase plane trajectory with leading order $\what(x,y,t)$ effects included. The red circle denotes the trajectory starting point and the green square denotes the trajectory end point after one period. Position of the trajectory end point $\rbar(T)$ as a function of (b) $\log_{10}({\rm Err})$ and (c) $M=L/\Delta x$ defining convergence.}
\label{fig:resolution}
\end{figure}
\begin{figure}[!p]
\begin{center}
(a)\includegraphics[width=0.4\textwidth]{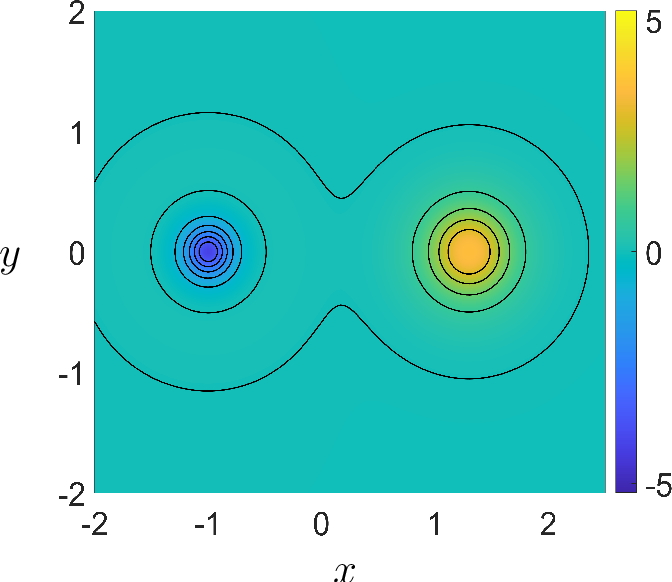}
(b)\includegraphics[width=0.4\textwidth]{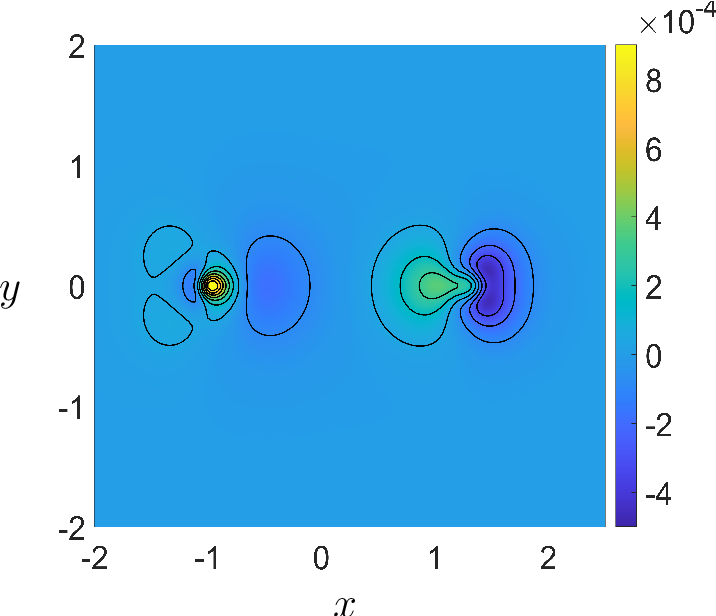}
(c)\includegraphics[width=0.4\textwidth]{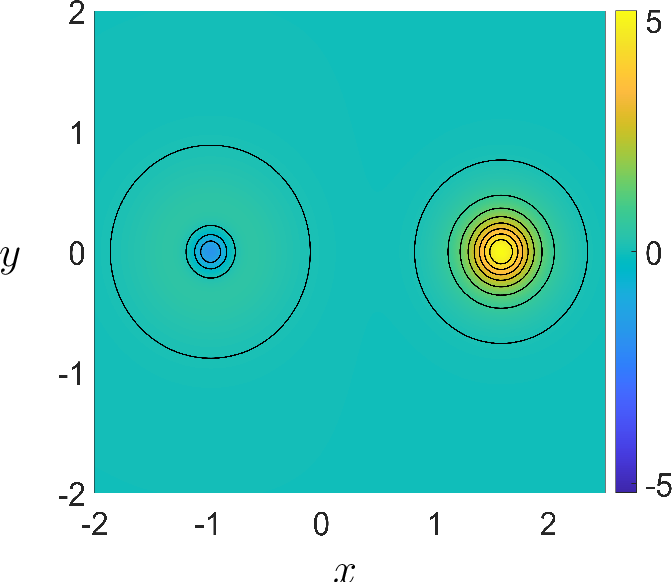}
(d)\includegraphics[width=0.4\textwidth]{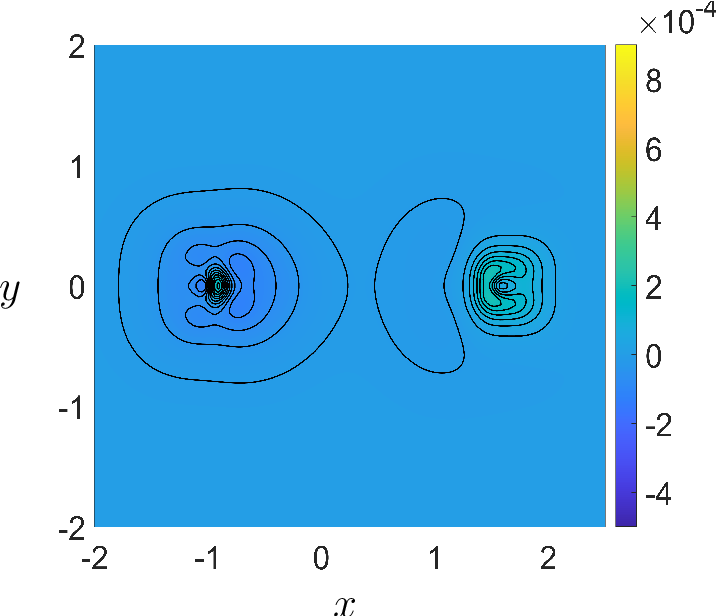}
(e)\includegraphics[width=0.4\textwidth]{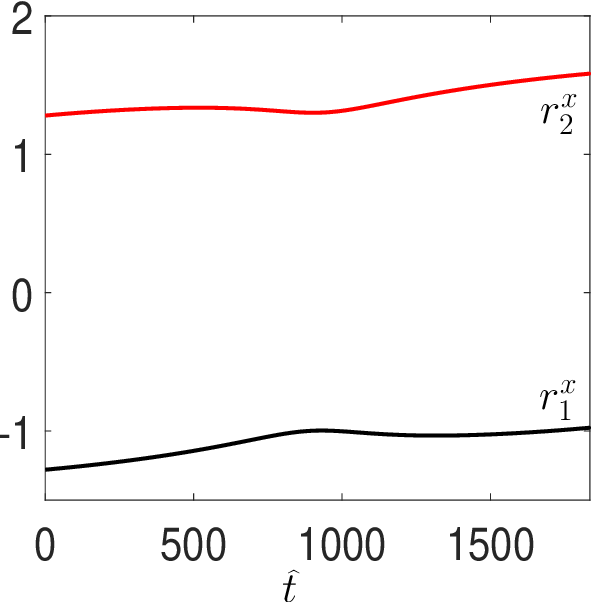}
\end{center}
\caption{Heat map of the real-part of the solution (a,c) $u(x,y)$ and the remainder function (b,d) $w(x,y)=\epsilon\what(x,y)$ for the trajectory starting with $(\rbar,\gbar)=(2.56,\pi/2)$ at (a,b) $\that=917.06$ where $\gbar=\pi/2$ and $\dfrac{\d\gbar}{\d\that}<0$ and (c,d) $\that=T=1834.96$ where $\gbar=\pi/2$ and $\dfrac{\d\gbar}{\d\that}>0$. (e) Plot of $r_1^x(\that)$ (black line) and $r_2^x(\that)$ (red line) for this trajectory showing the drift of the pulses to positive $x$ values over each period of the phase.}
\label{fig:what}
\end{figure}

When we consider the $O(\epsilon)$ effects of the projection scheme, i.e. when we retain the leading order effects of the remainder function $w=\epsilon\what$ as in the PS, then this periodic structure breaks down. In Figure \ref{fig:resolution}(a), we consider the phase trajectory for a result with initial separation $\rbar=2.56$ and the phase difference $\gbar=\pi/2$. Here the red circle signifies the starting point of the trajectory when $\gbar(0)=\pi/2$ and the green square denotes the end point of the trajectory when $\gbar(T)=\pi/2$ again, where $T$ denotes the time value (when it exists) where this phase difference reoccurs. This event can be detected using \lstinline[style=Matlab-editor]{MATLAB}'s \lstinline[style=Matlab-editor]{ode15s} subroutine. Note, the value of $T$ is not unique, and in fact there are multiple times when this condition is satisfied, so we also need to stipulate a condition on the sign of $\dfrac{\d\gbar}{\d\that}$ to uniquely determine this value. For $\dfrac{\d\gbar}{\d\that}>0$ the trajectory corresponds to $\that=T$ and the value of $\gbar=\pi/2$ again at the initial point of the trajectory, while $\dfrac{\d\gbar}{\d\that}<0$ with $\gbar=\pi/2$ corresponds to the point on the trajectory directly below the centre in Figure \ref{fig:phase_space}. The trajectory in Figure \ref{fig:resolution}(a) shows that at $\that=T$ the value of  $\rbar(T)-\rbar(0)<0$, which means the two pulses move closer together in time. A plot of the pulse positions along the trajectory where $\gbar=\pi/2$ with $\dfrac{\d\gbar}{\d\that}<0$ and $\dfrac{\d\gbar}{\d\that}>0$, along with the corresponding forms of $w(x,y)=\epsilon\what(x,y)$ are given in Figure \ref{fig:what}. The Figure shows that $|w|\ll|V|$ in the weakly interacting limit, and also demonstrates that the pulses slowly drift in space as they interact with each other.

Before we go on to discuss the significance of the pulses moving together over time, we first consider the convergence of the numerical scheme in Figures \ref{fig:resolution}(b,c). Here we plot the value of $\rbar(T)$ for the result in Figure \ref{fig:resolution}(a) as a function of (b) $\log_{10}({\rm Err})$, the error of the \lstinline[style=Matlab-editor]{MATLAB} integration scheme, and (c) the spatial grid size $M$, where $\Delta x=L/M$. It is clear from the scales on the $y$-axis that the results very quickly converge as ${\rm Err}$ and $\Delta x$ are reduced. In the results which we present in this paper, we use ${\rm Err}=10^{-8}$ and $M=150$ ($\Delta x=4\times 10^{-2}$), which give a good balance between resolution and run time.

\begin{figure}[!htb]
\begin{center}
(a)\includegraphics[width=0.45\textwidth]{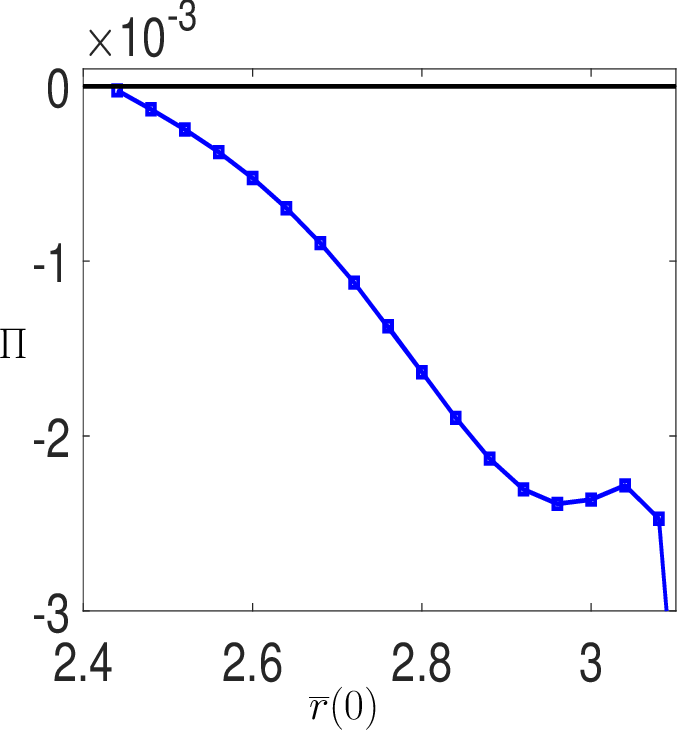}
(b)\includegraphics[width=0.45\textwidth]{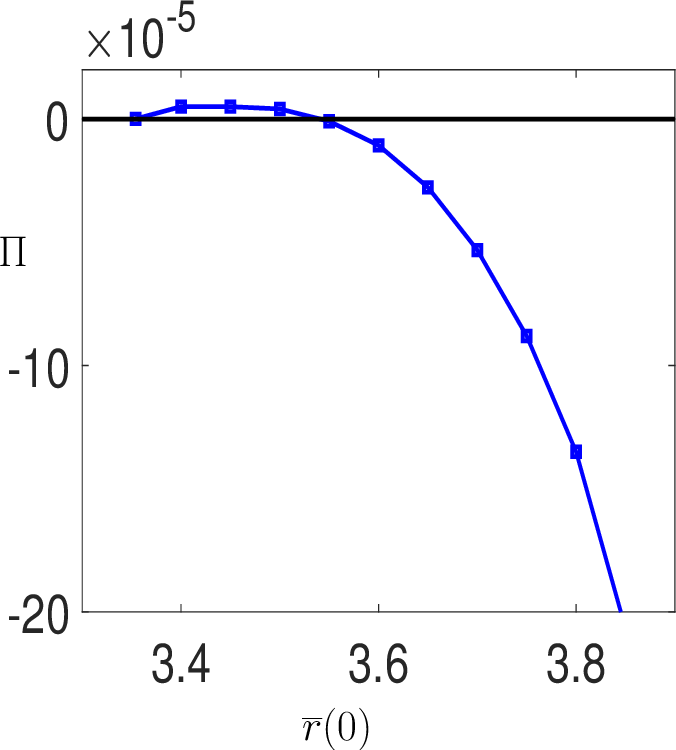}
\end{center}
\caption{Plot of $\Pi(\rbar(0))$ for the parameters in (a) (\ref{eqn:parameters1}) and (b) (\ref{eqn:parameters2}). In (b) the value of $\rbar(0)$ where $\Pi=0$ signifies the existence of a periodic orbit.}
\label{fig:periodic_orbit}
\end{figure}
By considering multiple trajectories around the equilibrium centre/focus point, like the one in Figure \ref{fig:resolution}(a), we can determine the dynamics of the system by considering the function
\begin{equation}
\Pi(\rbar(0))=\rbar(T)-\rbar(0).
\end{equation}
Ultimately, the system dynamics can be classified into three possible scenarios:
\begin{enumerate}
\item $\Pi(\rbar(0))>0$: Unstable dynamics, pulses move apart and away from the equilibrium point,
\item $\Pi(\rbar(0))<0$: Stable dynamics, pulses move together and toward the equilibrium point,
\item $\Pi(\rbar(0))=0$: Neutral dynamics, pulses oscillate in position in a periodic limit cycle.
\end{enumerate}
Here $\rbar(0)$ is chosen to lie above the centre/focus equilibrium point, which in the phase plane plots in Figure \ref{fig:phase_space} lie at (a) $(\rbar,\gbar)=(2.432,\pi/2)$ and (b) $(3.352,\pi/2)$ for the parameters in (\ref{eqn:parameters1}) and (\ref{eqn:parameters2}) respectively. In Figure \ref{fig:periodic_orbit}, we plot the function $\Pi(\rbar(0))$ for the parameters in (a) (\ref{eqn:parameters1}) and (b) (\ref{eqn:parameters2}). What this Figure shows is that for $\beta_r=-0.05$ in (a), $\Pi<0$ for all $\rbar(0)$ in cell 1, this shows that here the point at $(2.432,\pi/2)$ is the only stable point, and all trajectories will spiral into this stable focus. For $\beta_r=-2$ in (b) however, between the equilibrium point $(3.352,\pi/2)$ and the point $(3.543,\pi/2)$, the value of $\Pi>0$ and so the equilibrium point is an unstable focus and trajectories spiral out towards the limit cycle which exists at $\rbar(0)\approx 3.543$. This same qualitative behaviour, i.e. having a stable equilibrium transitioning to a stable limit cycle as $\beta_r$ varies, was also seen in the one-dimensional version of this problem, except here the values of $\beta_r$ where this transition occurs, is different \cite{rossides2023}. This quantitative difference is to be expected given the variations in the pulse shape seen in Figure \ref{fig:pulse}.

\subsection{Case $N>2$ pulse interactions}
\label{sec:N2greater}

\subsubsection{Steady states and basins of attraction}
\label{sec:basins}

For the case of $N>2$ pulses with no remainder function, the system no longer appears to be Hamiltonian as the governing equations cannot be determined from an energy functional, and hence the only way to identify information about the system is to directly integrate the governing system equations and seek states, such as steady fixed points or periodic states.  In cases where only translational pulse states are considered, i.e. all the phases are constant, then the system could be written in gradient form, as was the case in the two-dimensional  optical cavity problem of \cite{vladimirov2002}. However, we are unable to fix the phases in this way to make this simplification and so have to resort to numerical solutions of the system only.

\begin{figure}[!htb]
\begin{center}
\includegraphics[width=0.6\textwidth]{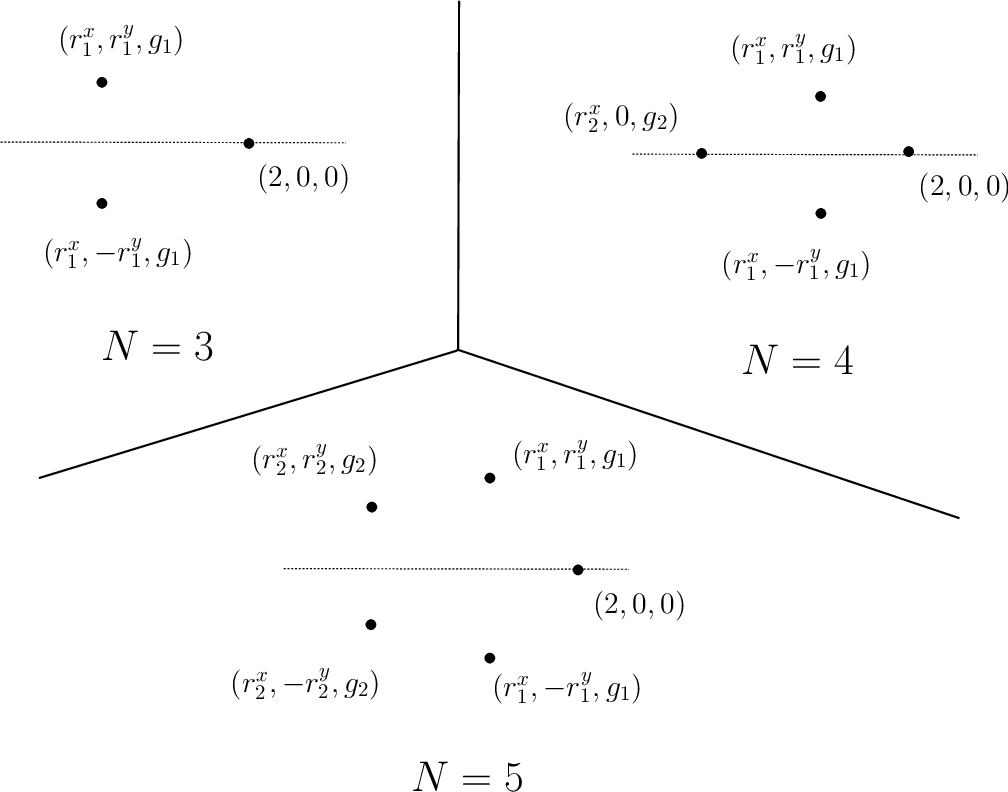}
\end{center}
\caption{Initial pulse configurations considered when seeking fixed point equilibria for $N=3,~4,~5$ pulses. The dotted horizontal lines denote the $x$-axis in each case.}
\label{fig:initial_fixed_config}
\end{figure}
Here we examine the properties of systems of $N=3,~4,$ and $5$ pulses. The first key piece of information to identify is whether there exists any stationary equilibria points in parameter space, which are solutions $(\rb_1,...,\rb_N)^T$ of (\ref{eqn:matrix}) with $\bigdot{{\bf X}}={\bf 0}$. In principle this is a challenging task numerically, because of the huge parameter space which needs to be explored. However, we can make progress by examining pulse configurations with particular symmetries as initial guesses, in particular we consider those regular geometric states given in the schematic Figure \ref{fig:initial_fixed_config}. In each case, we seek fixed equilibrium states relative to the other pulses. Hence, we are able to fix the steady state properties of one pulse, thus we choose to set
\begin{equation}
r_N^x=2,~~~~~r_N^y=0,~~~~~g_N=0,
\label{eqn:rN}
\end{equation}
and then seek solutions with reflectional symmetry through the $x$-axis (denoted by the dotted lines in Figure \ref{fig:initial_fixed_config}), in order to reduce the number of unknowns and to make the solution procedure more robust. Thus for $N=3$, have the three unknowns
\[
(r_1^x,r_1^y,g_1)~~~~~{\rm with}~~~~~(r_2^x,r_2^y,g_2)=(r_1^x,-r_1^y,g_1),
\]
for $N=4$ we, have the five unknowns
\[
(r_1^x,r_1^y,g_1,r_2^x,g_2)~~~~~{\rm with}~~~~~(r_2^y,r_3^x,r_3^y,g_3)=(0,r_1^x,-r_1^y,g_1),
\]
and for $N=5$ we, have the six unknowns
\[
(r_1^x,r_1^y,g_1,r_2^x,r_2^y,g_2)~~~~~{\rm with}~~~~~(r_3^x,r_3^y,g_3,r_4^x,r_4^y,g_4)=(r_2^x,-r_2^y,g_2,r_1^x,-r_1^y,g_1).
\]
In order to readily sweep the three-, five- and six-dimensional parameter spaces, we first identify fixed points by solving the $O(1)$ POS problem (\ref{eqn:O1}) with $\what\equiv0$, using Newton iterations. These fixed points are then used as initial guesses in the PS (\ref{eqn:Ow}) with $\what\neq0$ to identify fixed points of the full nonlinear system. Here we present results only for the parameter set given in (\ref{eqn:parameters1}). 

\begin{figure}[!htb]
\begin{center}
(a)\includegraphics[width=0.4\textwidth]{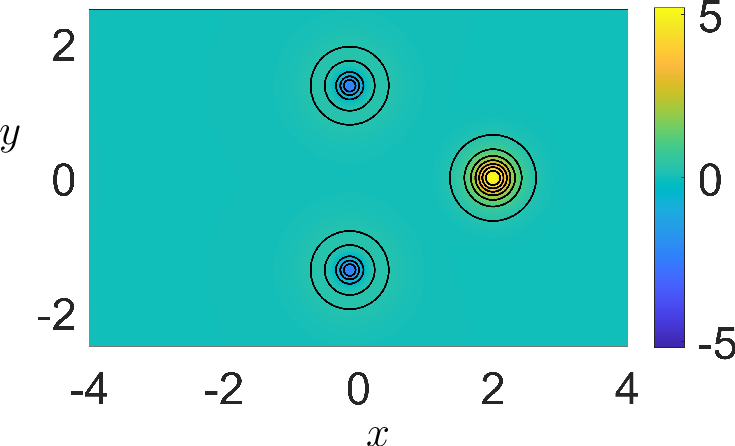}
(b)\includegraphics[width=0.4\textwidth]{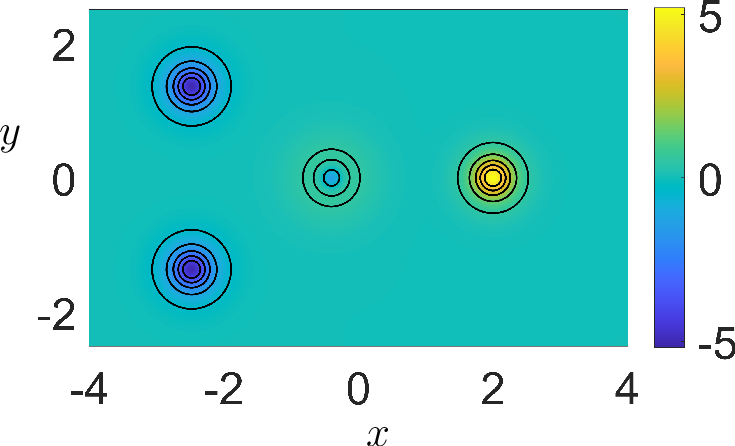}
(c)\includegraphics[width=0.4\textwidth]{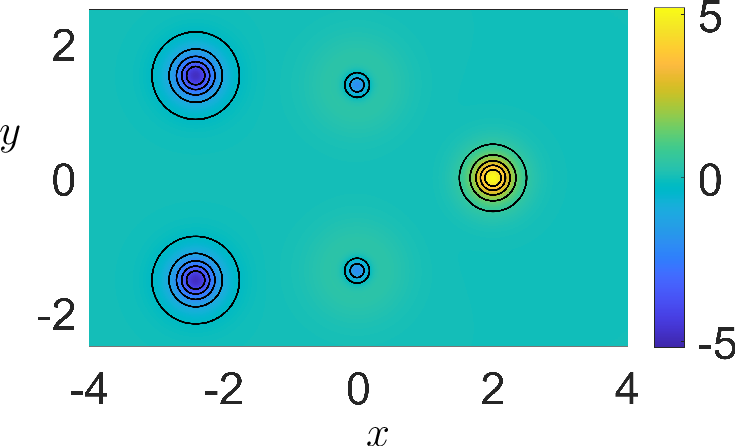}
\end{center}
\caption{Plot of ${\rm Re}(u)$ for the equilibrium point solutions for (a) $N=3$, (b) $N=4$ and (c) $N=5$ for the parameters in (\ref{eqn:parameters1})}
\label{fig:fixed_points}
\end{figure}

One set of possible equilibrium positions are the two-dimensional set, of pulses lined up along the $y$-axis with
\begin{eqnarray}
(r_1^x,r_1^y,g_1,J_{\rm max})&=&(2.00,1.38,~3.08,8.02\times10^{-1}), \nonumber\\
&=&(2.00,2.39,-1.76,2.49\times10^{-3}),\label{eqn:linedup}\\
&=&(2.00,3.27,-0.06,2.73\times10^{-6}).\nonumber
\end{eqnarray}
The term $J_{\rm max}$ denotes the largest eigenvalue of the corresponding Jacobian matrix of the Newton solver, and as $J_{\rm max}>0$ in each case, then these points are unstable to symmetric perturbations. Of interest to us in this study is the existence of stable equilibrium points $(J_{\rm max}<0)$ as these would constitute structures which persist for large times, under symmetric perturbations, and hence it would be a convergence point for solutions to (\ref{eqn:Ow}). In each of the cases, $N=3,~4,~5$, we could only identify one stable equilibrium point in the full $(\what\neq0)$ system. These stable point values are 
\begin{eqnarray}
&N=3&~~~~(r_1^x,r_1^y,g_1,J_{\rm max})=(-0.13,1.37,1.72,-8.38\times10^{-4}),\label{eqn:fixed3}\\
&N=4&~~~~(r_1^x,r_1^y,g_1,r_2^x,g_2,J_{\rm max})=(-2.48,1.36,-3.06,-0.40,1.48,-8.49\times10^{-4}),\label{eqn:fixed4}\\
&N=5&~~~~(r_1^x,r_1^y,g_1,r_2^x,r_2^y,g_2,J_{\rm max})=(-0.02,1.38,1.62,-2.42,,1.52,-2.98,-4.17\times10^{-4}),\label{eqn:fixed5}
\end{eqnarray}
and the form of real part of each solution, ${\rm Re}(u)$ are given in Figure \ref{fig:fixed_points}.
 These results show that the $N=3$ and $N=5$ equilibrium points form a regular polygon configuration, while the $N=4$ result resembles a triangle of pulses with an additional  central pulse. Note, we were unable to identify stable equilibrium points for $N>5$ for initial guess configurations of the form similar to those in Figure \ref{fig:initial_fixed_config}, however it is possible that stable points may exist for more elaborate initial guess configurations.

\begin{figure}[!htb]
\begin{center}
(a)\includegraphics[width=0.4\textwidth]{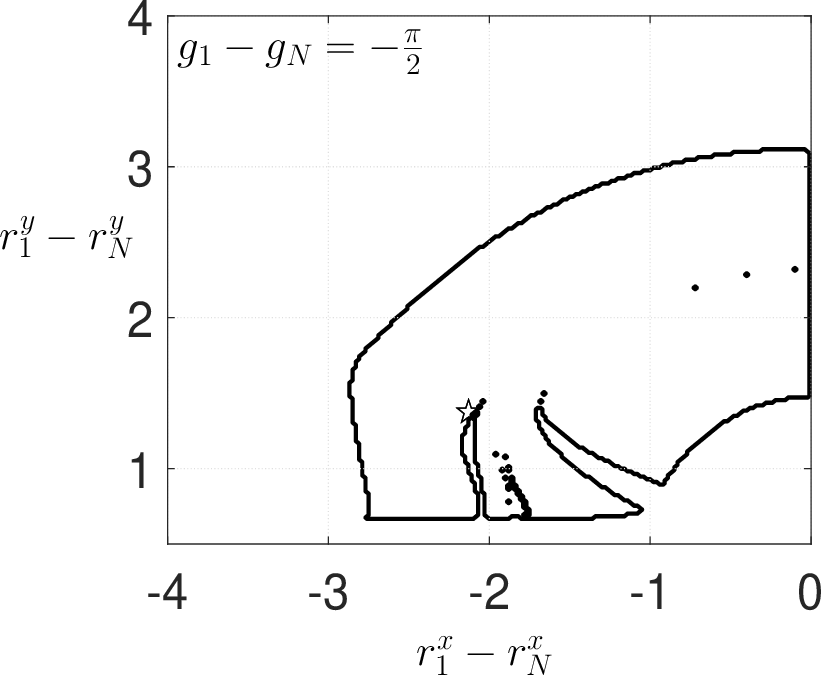}
(b)\includegraphics[width=0.4\textwidth]{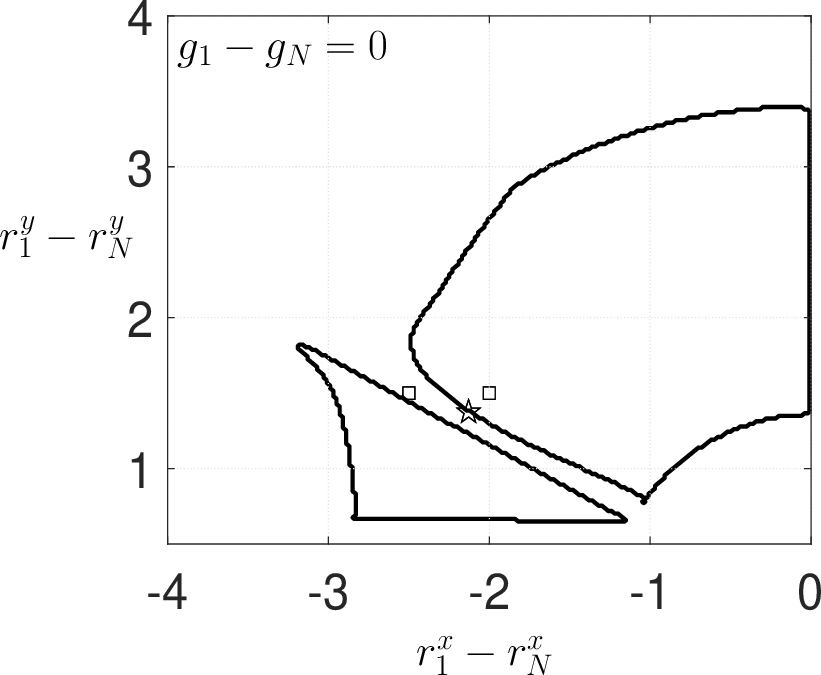}
(c)\includegraphics[width=0.4\textwidth]{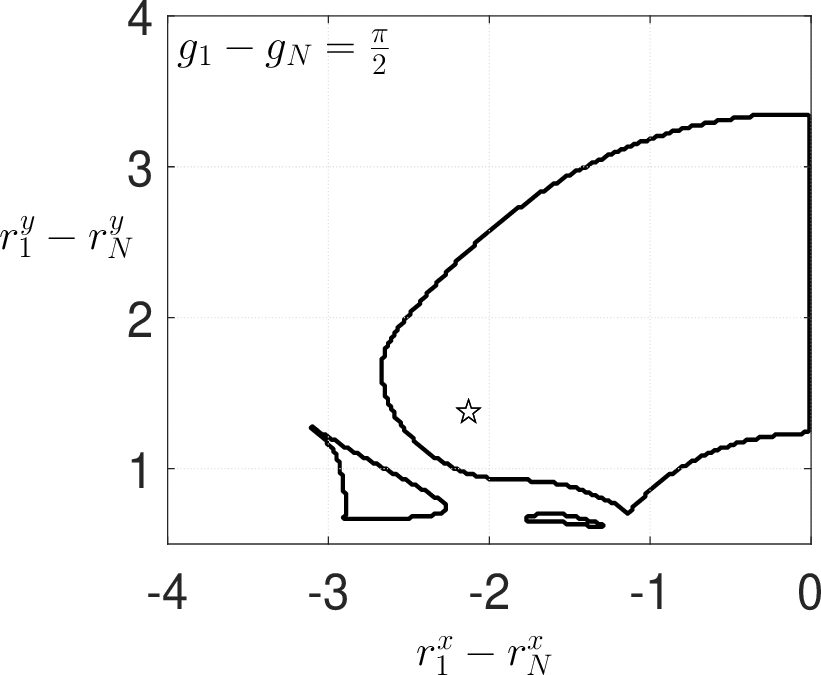}
(d)\includegraphics[width=0.4\textwidth]{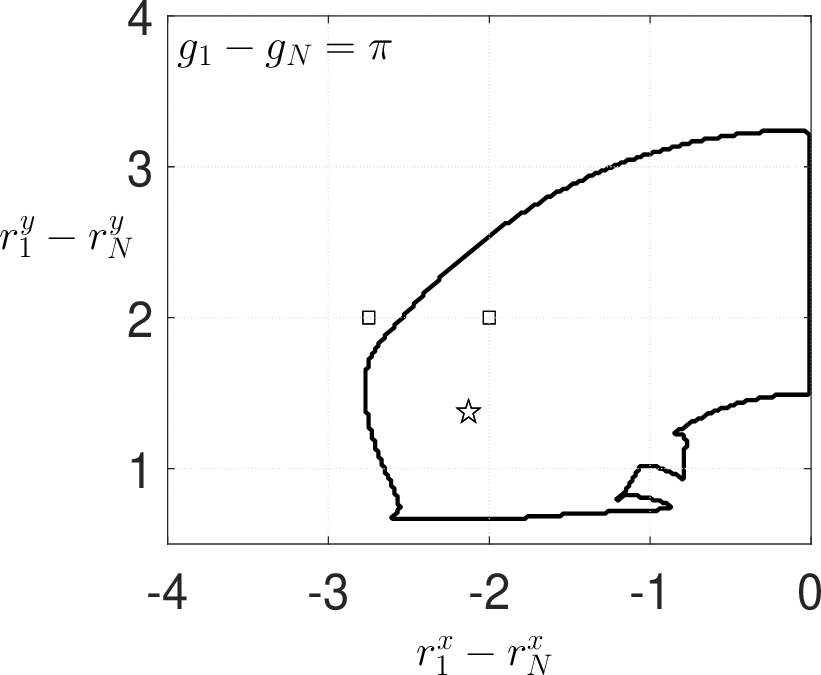}
\end{center}
\caption{Basin of attraction (inside black contour) for the $N=3$ fixed point as a function of $r_1^x-r_N^x$ and $r_1^y-r_N^y$ for (a) $g_1-g_N=-\pi/2$, (b) $g_1-g_N=0$, (c) $g_1-g_N=\pi/2$ and (d) $g_1-g_N=\pi$. The star denotes the projection of the equilibrium point (\ref{eqn:fixed3}) on to each plane, and the squares denote results considered in Figure \ref{fig:basin_N3_trajectories}. }
\label{fig:basin_N3}
\end{figure}
We expect the stable equilibrium points in (\ref{eqn:fixed3})-(\ref{eqn:fixed5}) above to be convergence points for the full system when some particular initial conditions are chosen. Here, we examine precisely the size and shape of the basin of attraction of these equilibrium points in parameter space. To calculate the basin of attraction, we solve the $O(1)$ POS version of the problem (\ref{eqn:O1}), as this amounts to solving ODEs only, which can be performed readily and to a high degree of accuracy (${\rm Err}=10^{-12}$). In order to reduce the number of initial condition parameters which need to be varied, we choose to fix the position of the $N^{\rm th}$ pulse to the values in (\ref{eqn:rN}), fix the other initial parameters to be their equilibrium point values from (\ref{eqn:fixed3})-(\ref{eqn:fixed5}), i.e. $(r_2^x,g_2)$ for $N=4$, $(r_2^x,r_2^y,g_2)$ for $N=5$, and hence only vary the three parameters $(r_1^x,r_1^y,g_1)$. Also, for both the POS and PS systems in (\ref{eqn:O1}) and (\ref{eqn:Ow}) respectively, we enforce the reflection symmetry given in Figure \ref{fig:initial_fixed_config}. Enforcing this symmetry means that when the solution reaches the equilibrium point, it remains there for a significant period of time. If we do not enforce this symmetry, then the pulses would only remain at the equilibrium point for a short period of time, because round off error in the calculation would ultimately build up in an asymmetric way, eventually causing the pulses to move away from the equilibrium point, and exhibit chaotic motions. This is because these equilibrium points are only stable to symmetric perturbations, by construction, not to asymmetric ones.

In Figure \ref{fig:basin_N3}, we plot the basin of attraction of the equilibrium point (\ref{eqn:fixed3}) in the $(r_1^x-r_N^x)(r_1^y-r_N^y)$-plane for different values of $g_1-g_N$. Here the ODES (\ref{eqn:O1}) are integrated to $\that=6\times10^{4}$. We check the values of $|r_1^x-r_N^x|,~|r_1^y-r_N^y|$ and $|g_1-g_N|$ at this point, and if they are all found to be within $10^{-4}$ of the fixed point values, then we declare the point to lie in the basin. The results presented here were checked against using different end times, and the results were consistent with each other. The size of the basin is large, and relatively coherent, although it appears to split into multiple regions for $g_1-g_N=0$ and $\pi/2$. However, in the three-dimensional parameter space these regions are all connected. For $g_1-g_N=-\pi/2$ we see some evidence of speckled behaviour, especially close to the basin boundary and this is evidence of the existence of chaos in the system at this point \cite{ott1994,mason2009,wright2015}. The basin extends close to $r_1^x-r_1^N=0$ where the only equilibrium points are those aligned pulses in (\ref{eqn:linedup}), which are unstable to symmetric perturbations.

\begin{figure}[!p]
\begin{center}
(a)\includegraphics[width=0.3\textwidth]{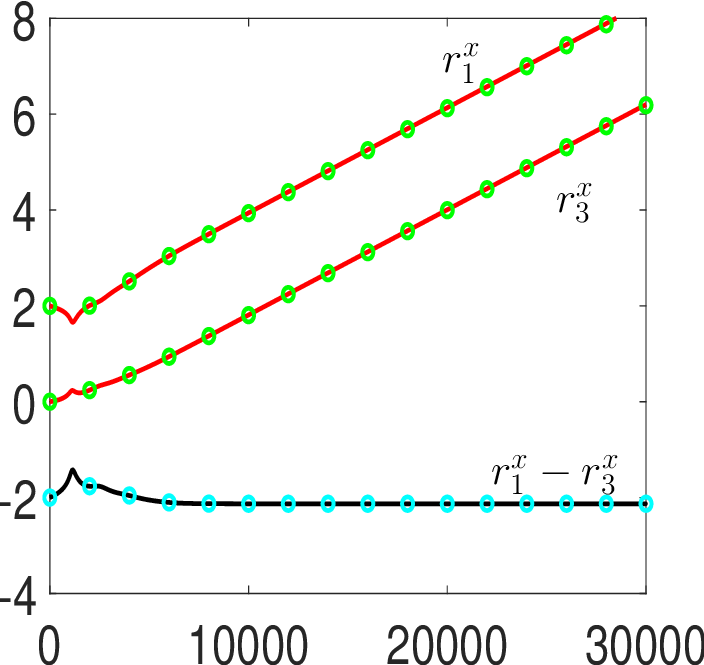}
\includegraphics[width=0.3\textwidth]{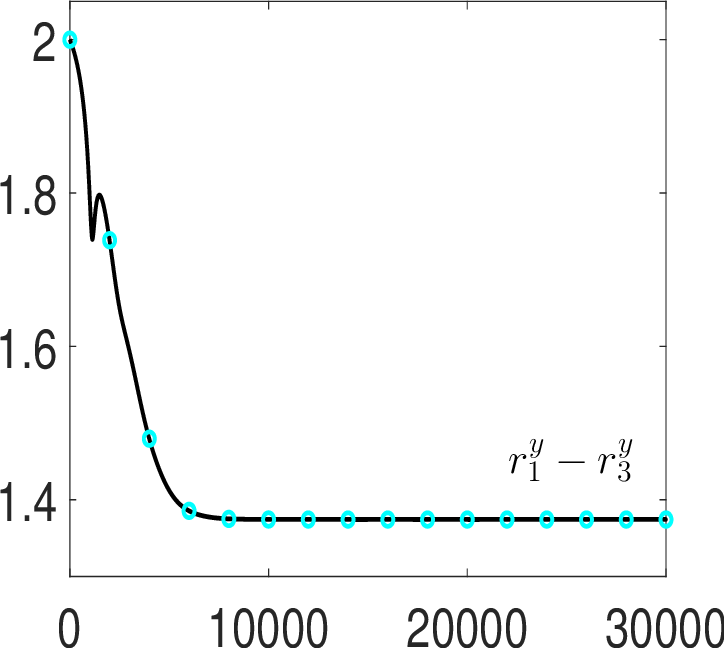}
\includegraphics[width=0.3\textwidth]{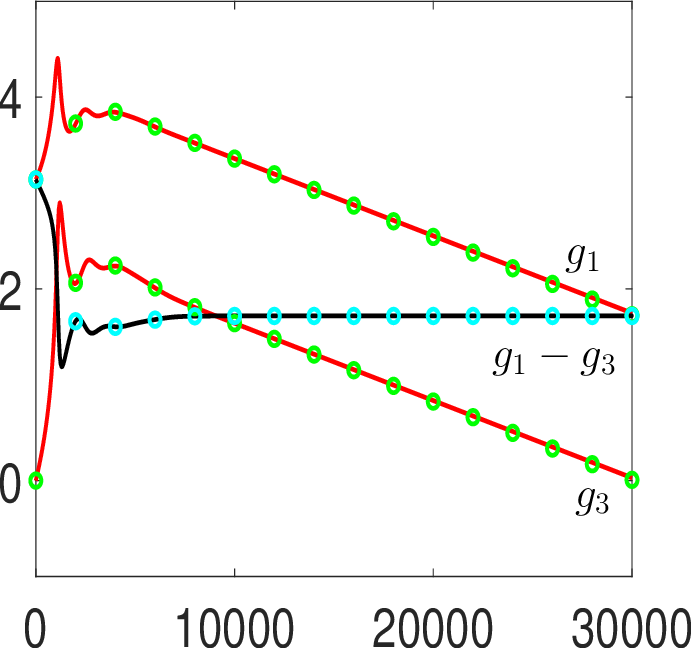}
(b)\includegraphics[width=0.3\textwidth]{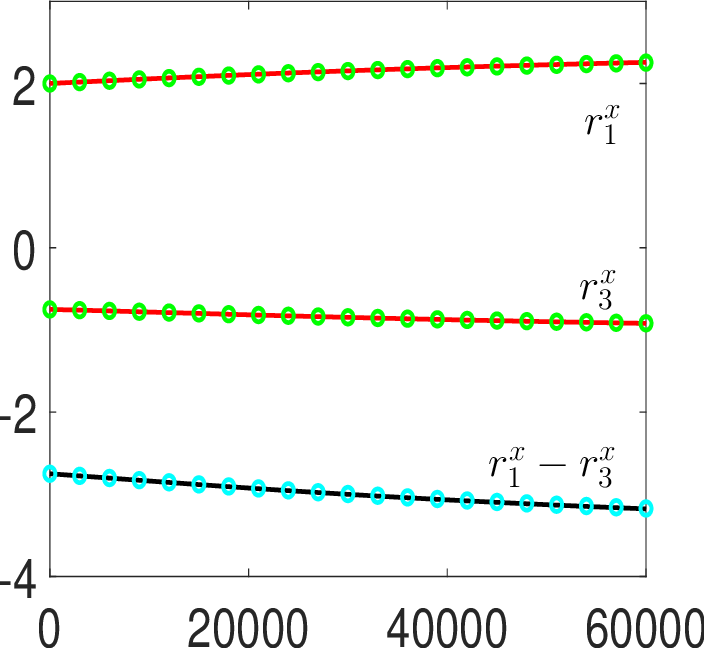}
\includegraphics[width=0.3\textwidth]{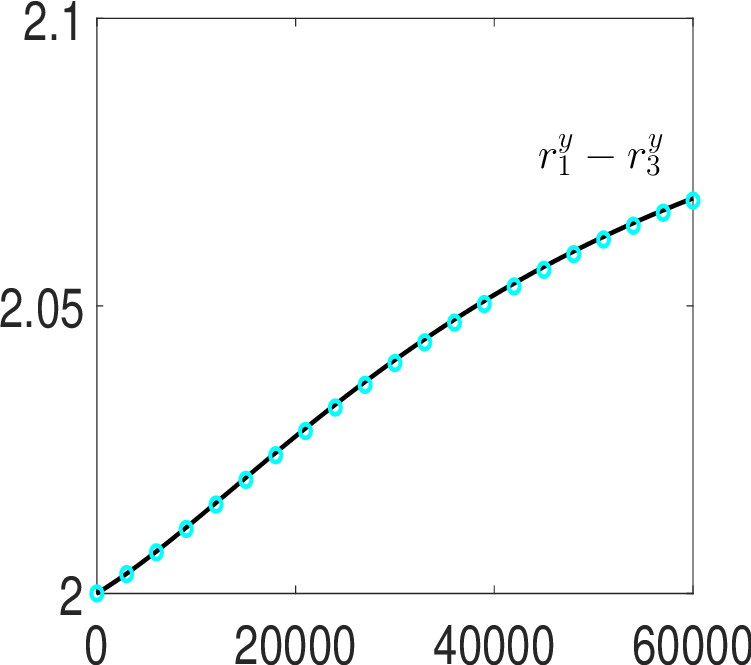}
\includegraphics[width=0.3\textwidth]{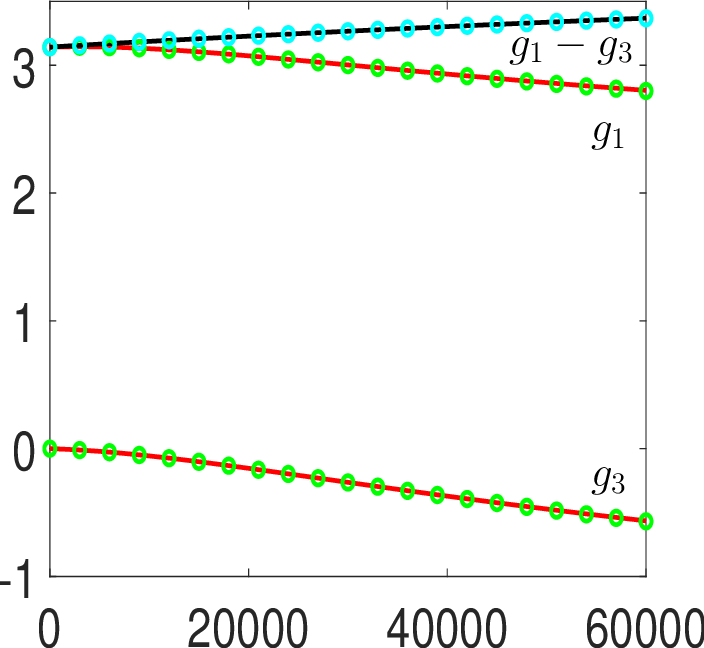}
(c)\includegraphics[width=0.3\textwidth]{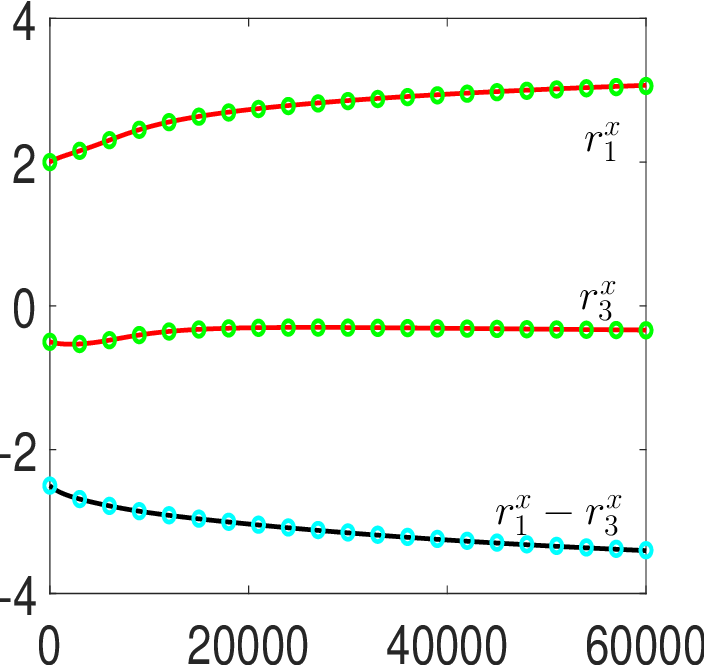}
\includegraphics[width=0.3\textwidth]{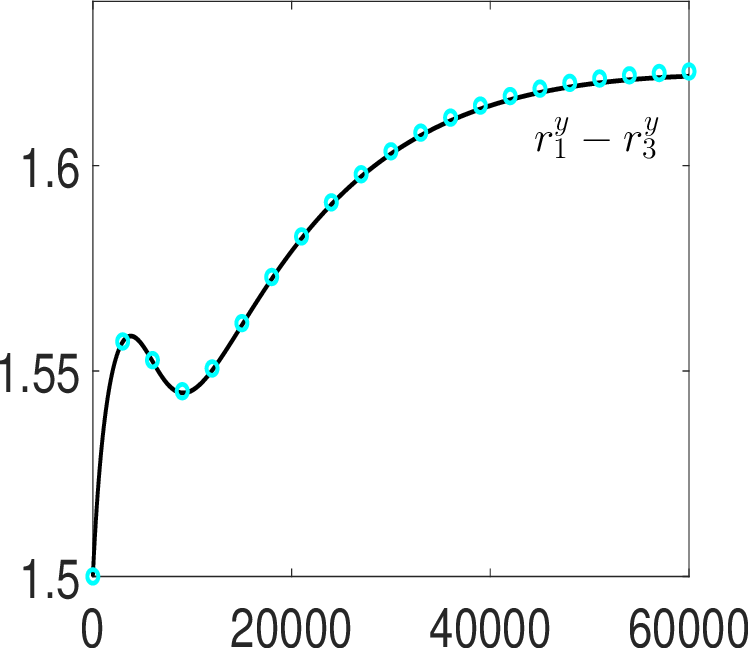}
\includegraphics[width=0.3\textwidth]{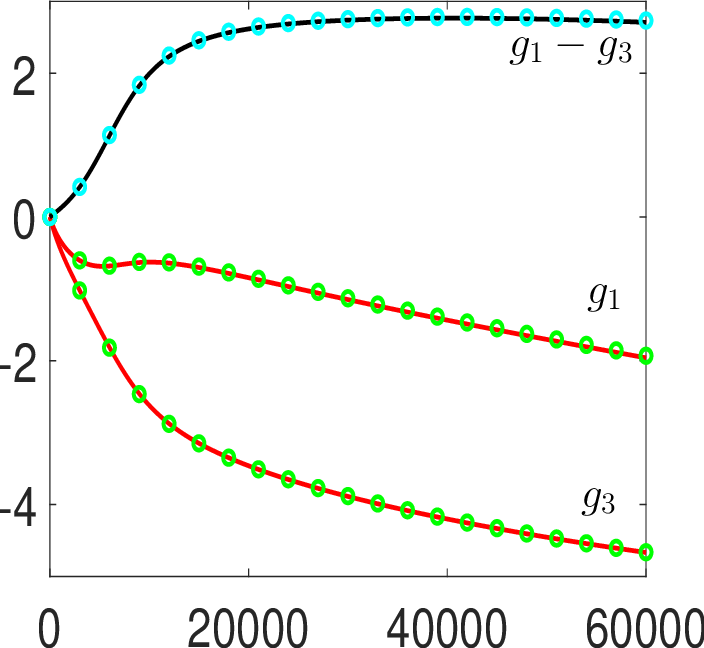}
(d)\includegraphics[width=0.3\textwidth]{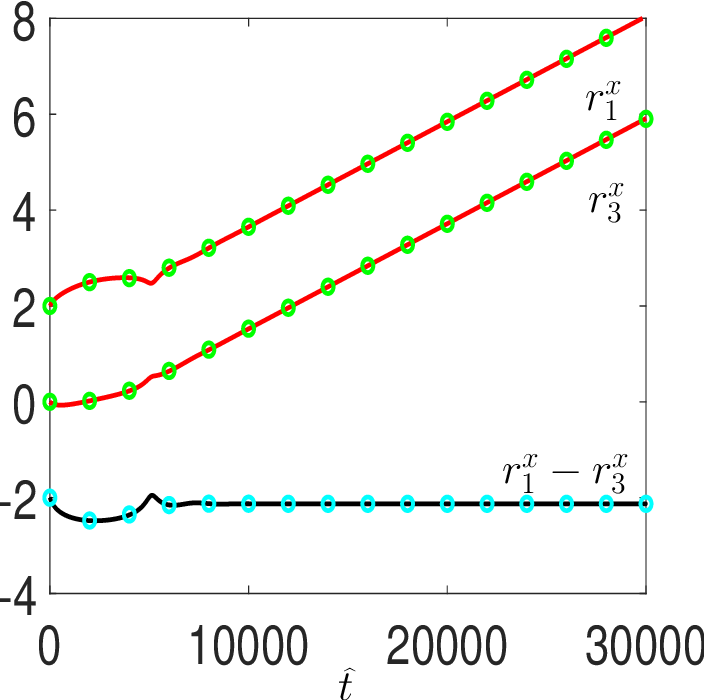}
\includegraphics[width=0.3\textwidth]{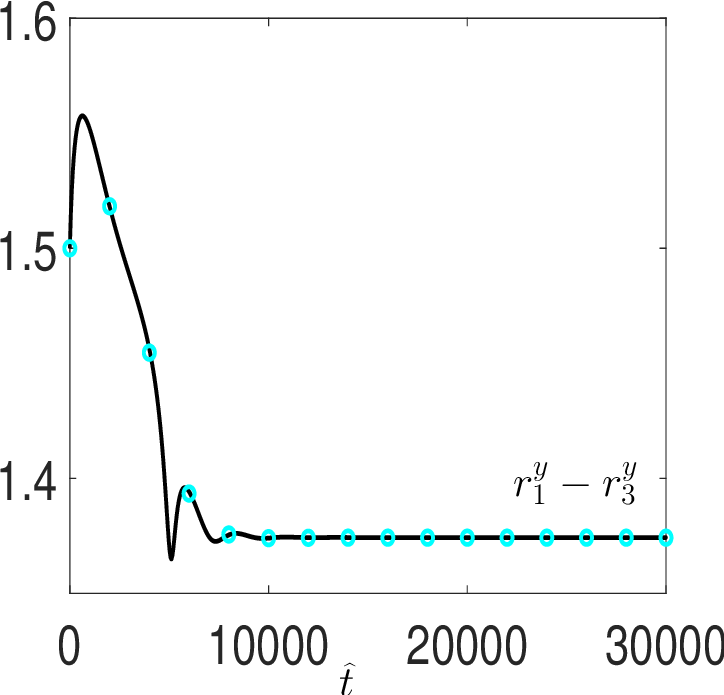}
\includegraphics[width=0.3\textwidth]{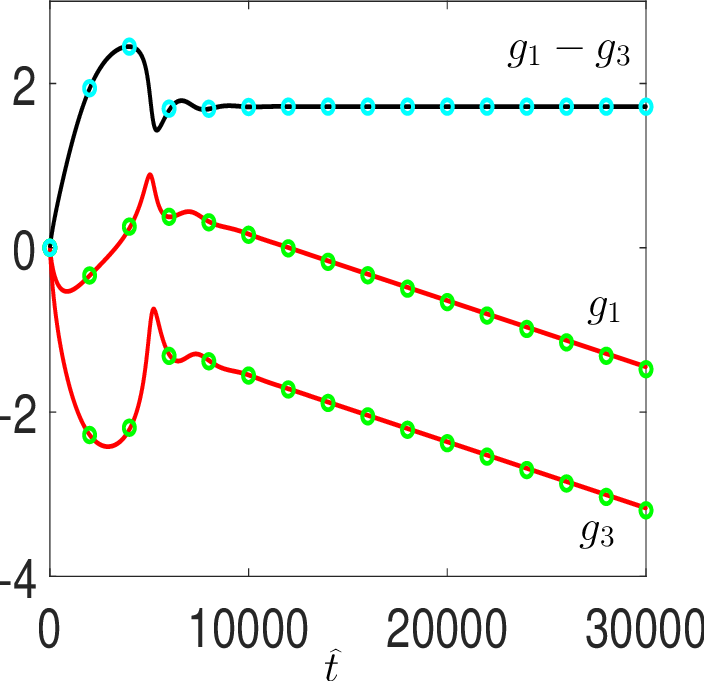}
\end{center}
\caption{Time-evolution plots of the pulse parameters for the case $N=3$ with the initial data (a) $(r_1^x,r_1^y,g_1)=(0,2,\pi)$, (b) $(-0.75,2,\pi)$, (c) $(-0.5,1.5,0)$ and (d) $(-0,1.5,0)$. In each case, the other parameters are given by (\ref{eqn:rN}). The solid lines are results of the POS (\ref{eqn:O1}) with $\what\equiv0$ and the symbols are the results of the PS (\ref{eqn:Ow}) with $\what\neq0$. The red lines signify the physical positions and phases $r_j^x,g_j$ of the pulses, while the black lines give the differences $r_1^x-r_3^x,r_1^y-r_3^y,g_1-g_3$. Fixed points are reached when these differences are constants.}
\label{fig:basin_N3_trajectories}
\end{figure}
The four squares in panels (b) and (d) indicate four cases, the time evolution of whose parameters are given in Figure \ref{fig:basin_N3_trajectories}. For each value of $g_1-g_N$ the two cases correspond to one inside and one outside the basin of attraction. For the cases inside the basin (Figures \ref{fig:basin_N3_trajectories}(a) and \ref{fig:basin_N3_trajectories}(d)) the three pulses quickly converge to the equilibrium point result, after an initial transient regime, and then propagate together in the positive $x$-direction. For the two cases outside the basin, the pulses appear to drift slowly away from one another, with the distance in the $y$-direction increasing most rapidly. These results have clearly not settled down to the equilibrium point result by $\that=6\times10^{4}$. Integrating these cases further we find that they do not converge to the equilibrium point in (\ref{eqn:fixed3}), nor to any other stable state, before numerical round-off begins contaminating the results. Hence, we are confident these cases lie outside the basin of attraction and that the results presented in Figure \ref{fig:basin_N3} are robust. The circles in Figure \ref{fig:basin_N3_trajectories} give the corresponding results for the time evolution of the pulses for the PS method (\ref{eqn:Ow}) when $\what\neq0$. For these results, the agreement with the leading order calculations is excellent; this is due to the fact that for these initial conditions the pulses do not leave the weakly-interacting regime.

\begin{figure}[!htb]
\begin{center}
(a)\includegraphics[width=0.4\textwidth]{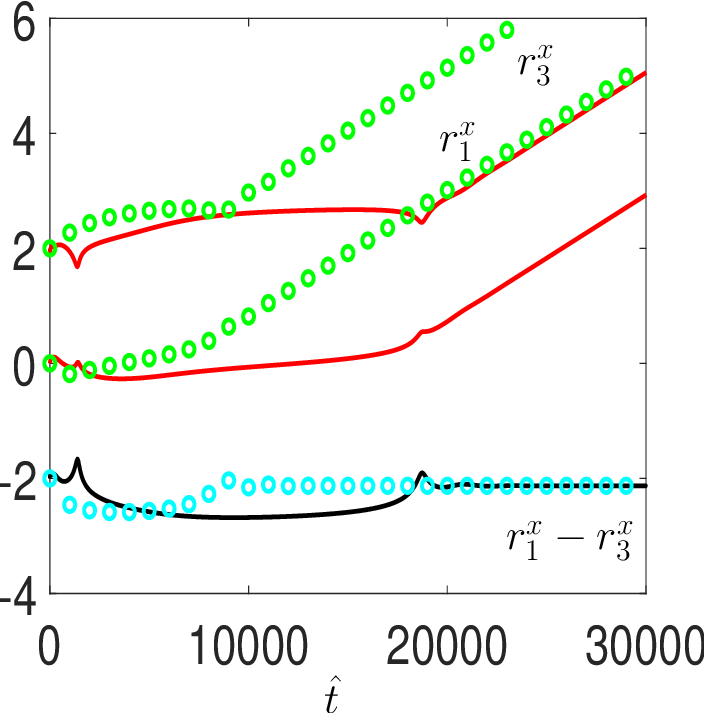}
\includegraphics[width=0.4\textwidth]{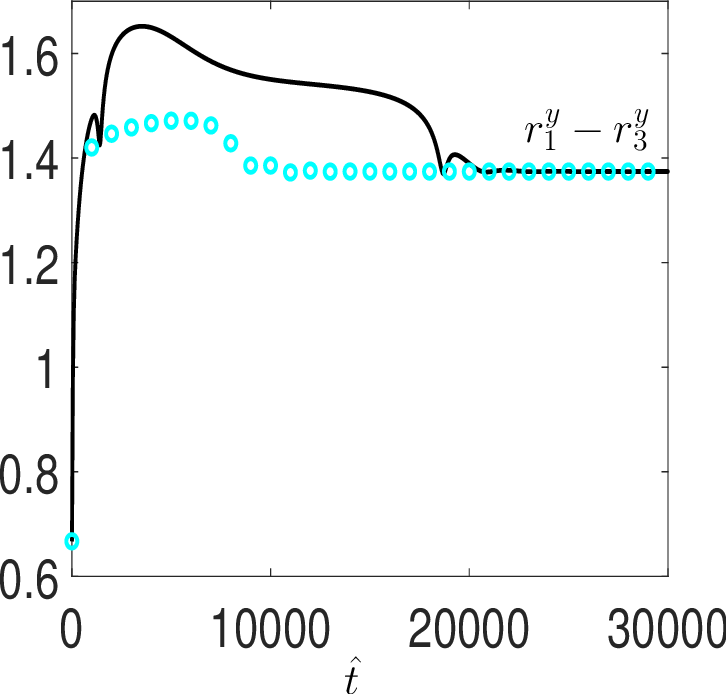}
(b)\includegraphics[width=0.4\textwidth]{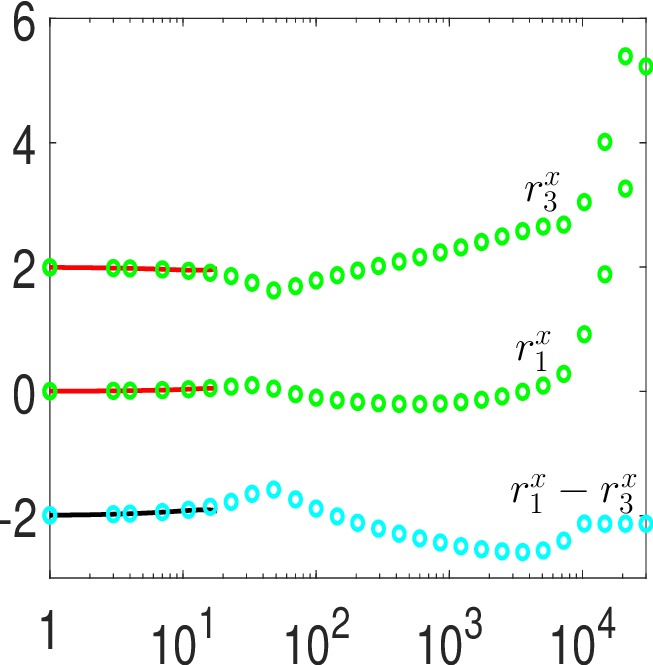}
\includegraphics[width=0.4\textwidth]{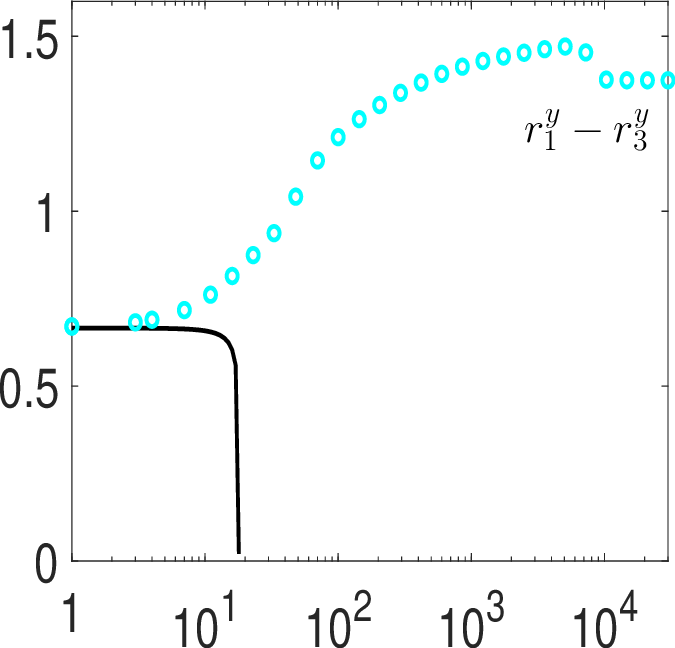}
\end{center}
\caption{Time-evolution plots of the pulse parameters for the case $N=3$ with the initial data (a) $(r_1^x,r_1^y,g_1)=(0,0.667,0)$ and (b) $(0,0.666,0)$. In each case, the other parameters are given by (\ref{eqn:rN}). The solid lines are results of the POS (\ref{eqn:O1}) with $\what\equiv0$ and the symbols are the results of the PS (\ref{eqn:Ow}) with $\what\neq0$. The red lines signify the physical positions $r_j^x$, while the black lines give the differences $r_1^x-r_3^x,r_1^y-r_3^y$.}
\label{fig:basin_N3_trajectories2}
\end{figure}
When we consider initial configurations where the pulses are close together (strongly interacting) then we can see differing results between the POS and PS results, similar to what we found for the $N=2$ pulse case in \S\ref{sec:N2}. In Figure \ref{fig:basin_N3_trajectories2}, we consider the time evolution for two initial configurations close to the basin boundary, with $(r_1^x,r_1^y,g_1)=(0,0.667,0)$ in (a) and $(0,0.666,0)$ in (b). The case in Figure \ref{fig:basin_N3_trajectories2}(a) lies inside the basin of attraction, and both the POS and PS methods integrate to the equilibrium point solution, but the PS scheme has a much shorter transient period, and a much smaller overshoot in $r_1^y-r_3^y$ compared to the leading order POS scheme. For the case in Figure \ref{fig:basin_N3_trajectories2}(b), the POS sees $(r_1^y-r_3^y)\to0$ very rapidly around $\that=10$, and  the two pulses annihilate. The PS scheme on the other hand converges to the equilibrium point in much the same way as the $r_1^x=0.667$ result in Figure \ref{fig:basin_N3_trajectories2}(b) (although with $\that$ on a log scale this is difficult to see), and hence this shows that the true basin of attraction for the full system is larger than is presented here, with the extension occurring close to regions where the pulses initially interact more strongly.

\begin{figure}[!htb]
\begin{center}
(a)\includegraphics[width=0.4\textwidth]{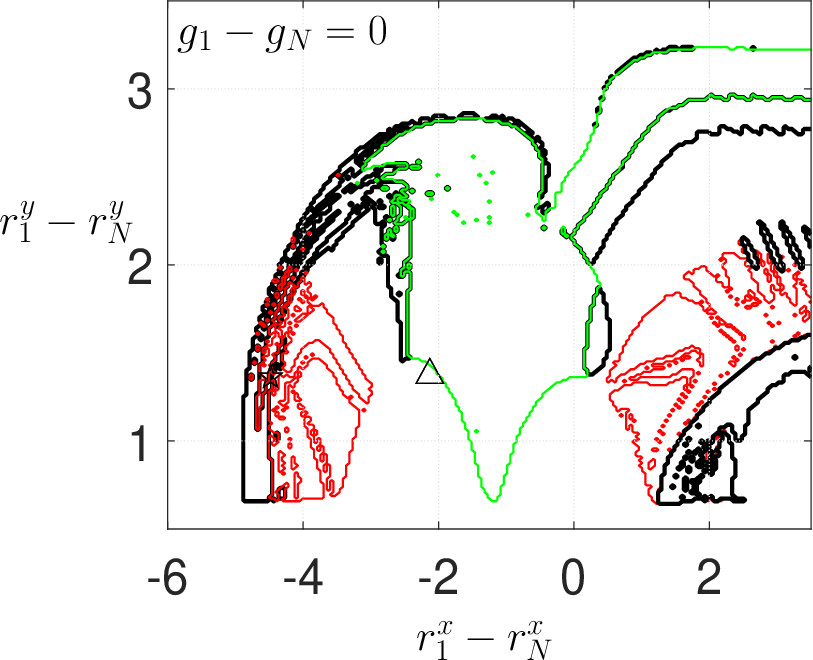}
(b)\includegraphics[width=0.4\textwidth]{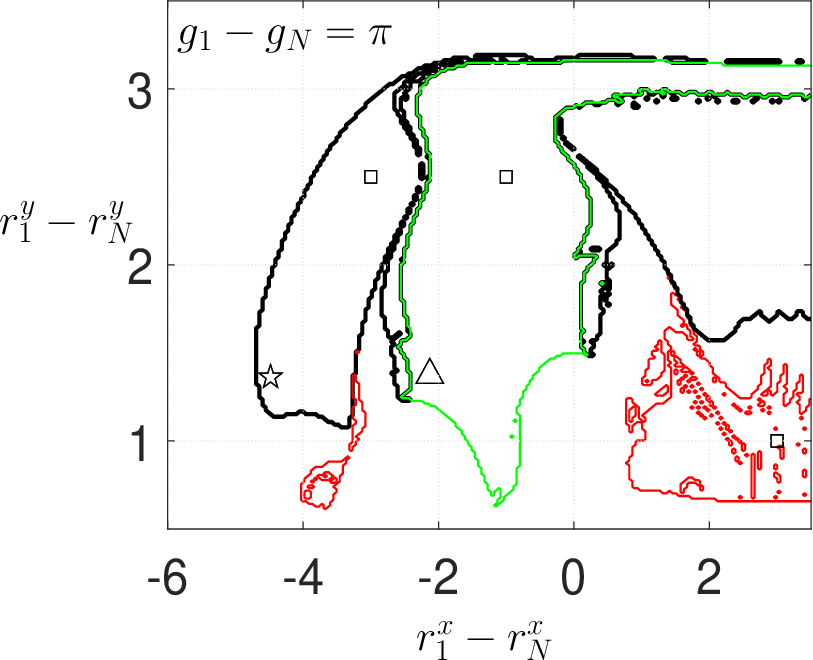}
(c)\includegraphics[width=0.4\textwidth]{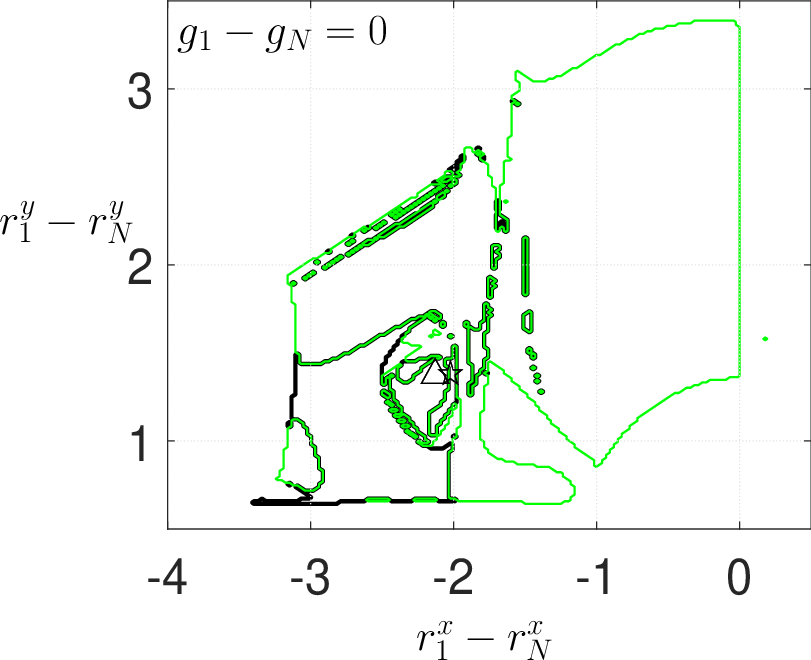}
(d)\includegraphics[width=0.4\textwidth]{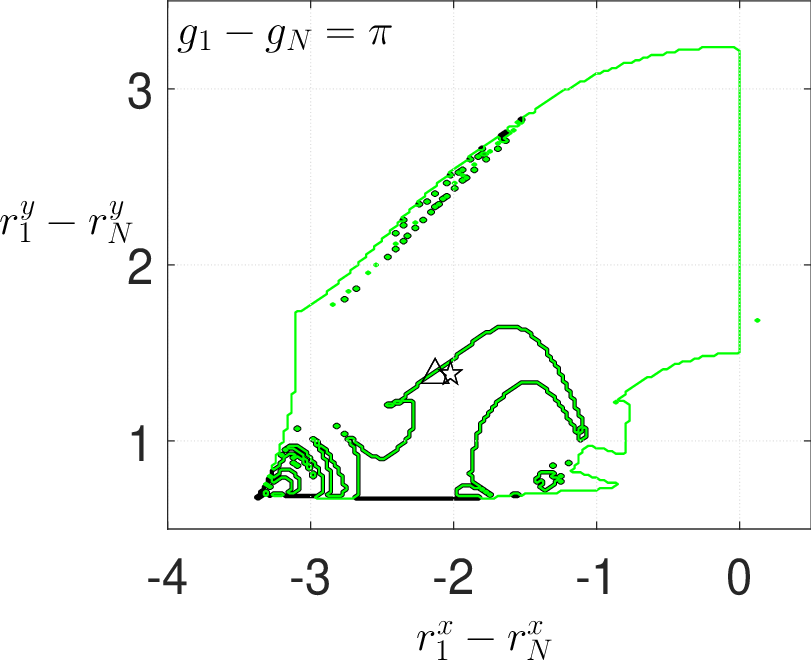}
\end{center}
\caption{Basins of attraction for (a,b) $N=4$ and (c,d) $N=5$ fixed point as a function of $r_1^x-r_N^x$ and $r_1^y-r_N^y$ for (a,c) $g_1-g_N=0$ and (b,d) $g_1-g_N=\pi$. The other initial parameters are given by the fixed point values in (\ref{eqn:rN}) together with (\ref{eqn:fixed4}) and (\ref{eqn:fixed5}) respectively. Inside the black contours the initial conditions converge to the $N=4,5$ equilibrium point respectively, inside the green contour they converge to the $N=3$ equilibrium point, and inside the red contour they converge to a periodic limit cycle. The star denotes the projection of the fixed point (\ref{eqn:fixed4}) and (\ref{eqn:fixed5}) respectively, while the triangle denotes the $N=3$ fixed point projection. The three squares in (b) denote results considered in Figure \ref{fig:basin_N4_trajectories}. }
\label{fig:basin_N45}
\end{figure}
\begin{figure}[!htb]
\begin{center}
(a)\includegraphics[width=0.3\textwidth]{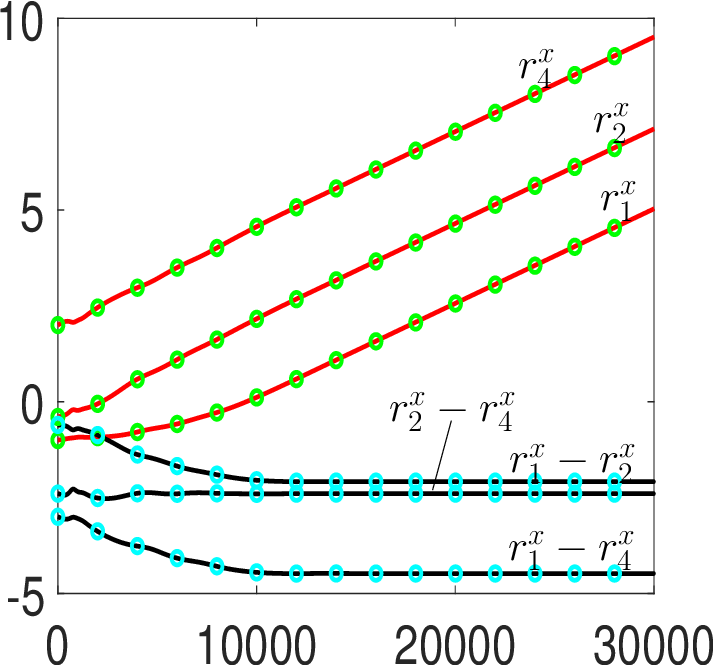}
\includegraphics[width=0.3\textwidth]{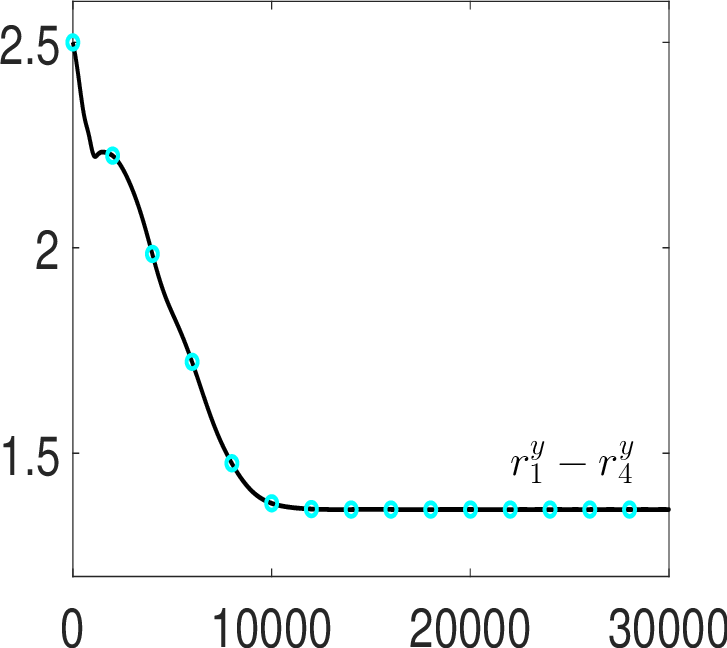}
\includegraphics[width=0.3\textwidth]{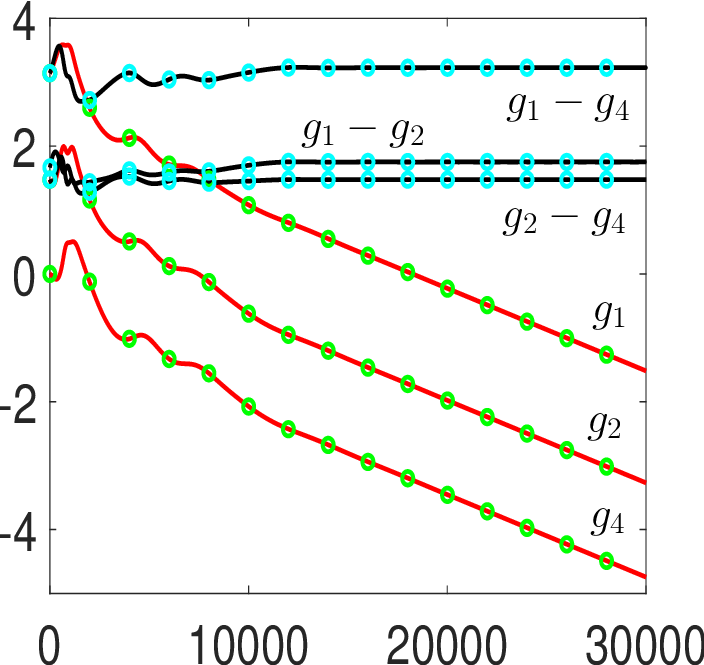}
(b)\includegraphics[width=0.3\textwidth]{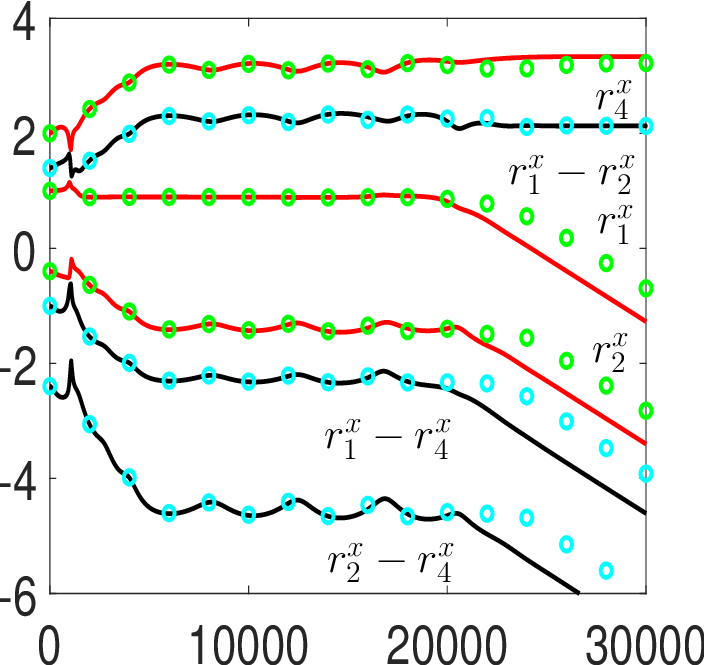}
\includegraphics[width=0.3\textwidth]{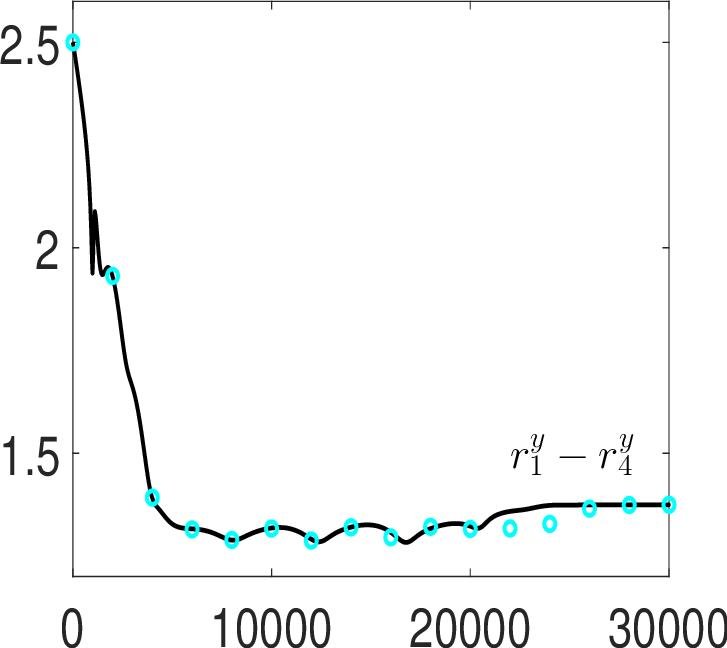}
\includegraphics[width=0.3\textwidth]{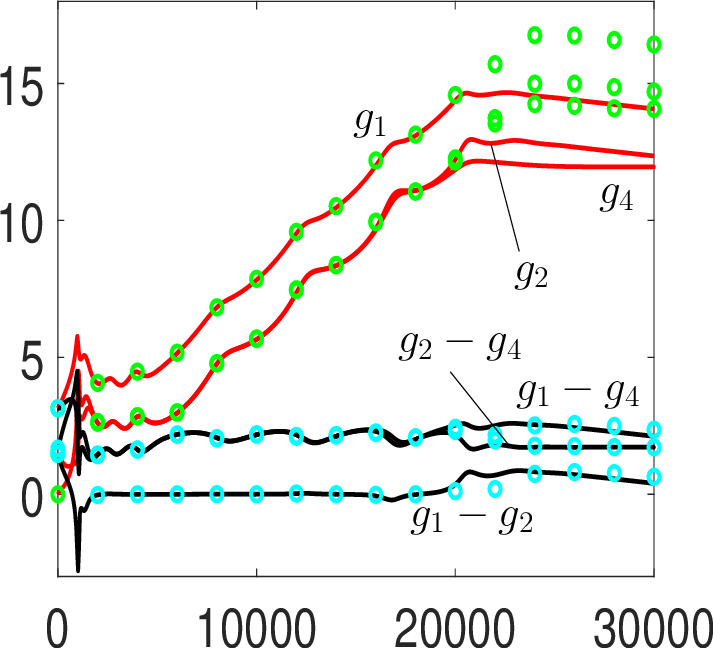}
(c)\includegraphics[width=0.3\textwidth]{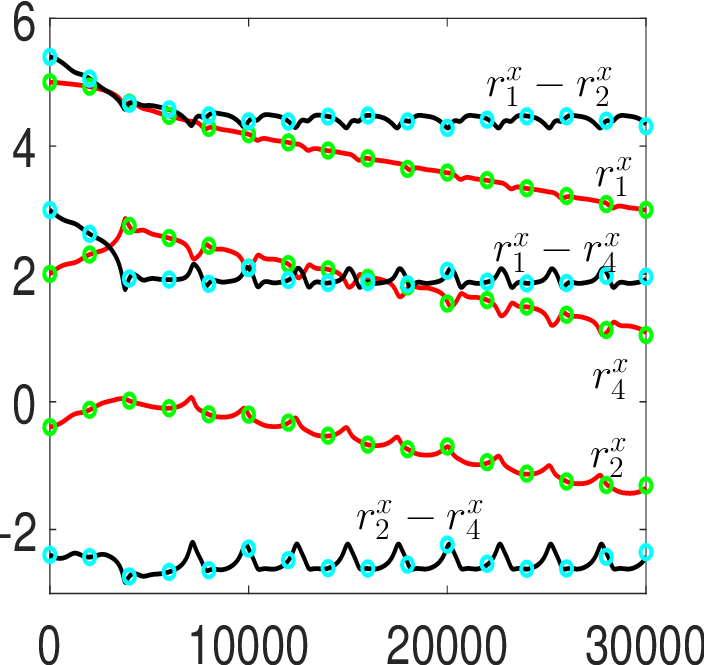}
\includegraphics[width=0.3\textwidth]{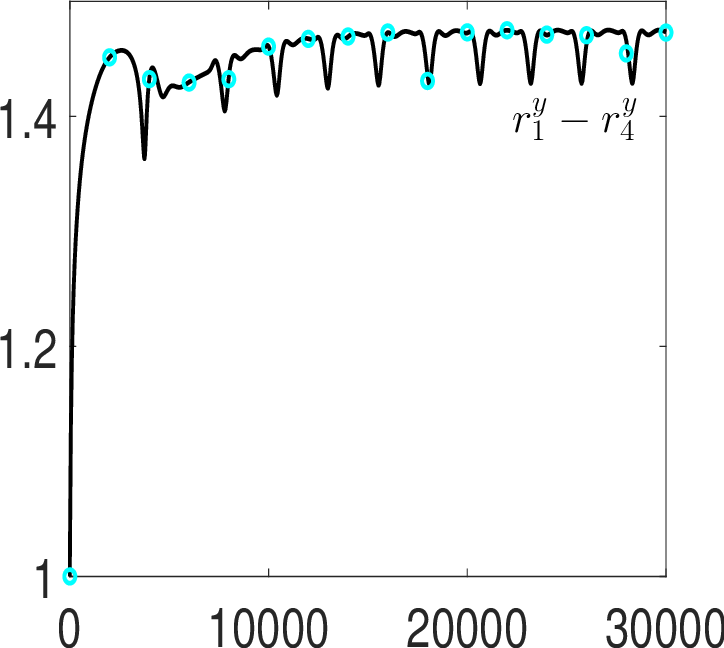}
\includegraphics[width=0.3\textwidth]{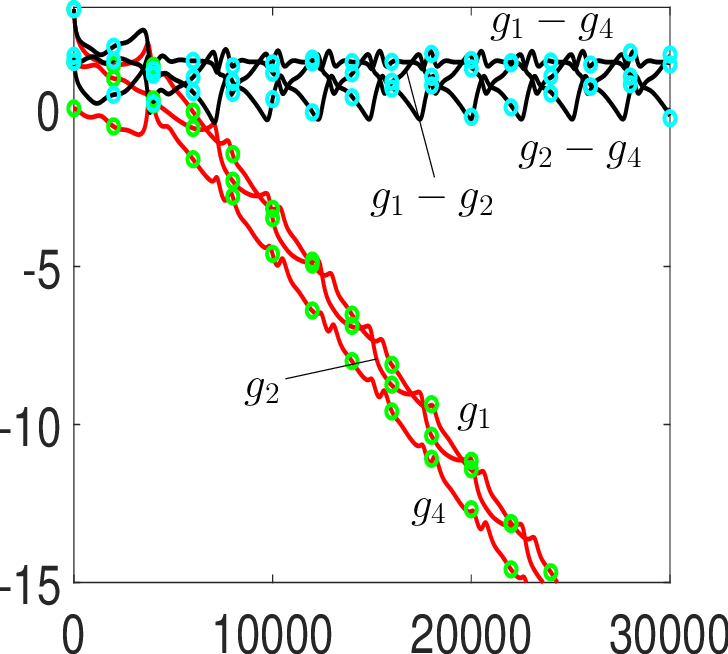}
\end{center}
\caption{Time-evolution plots of the pulse parameters for the case $N=4$ with the initial data (a) $(r_1^x,r_1^y,g_1)=(-1,2.5,\pi)$, (b) $(1,2.5,\pi)$ and (c) $(5,1,\pi)$. In each case the other parameters are given by the fixed point values (\ref{eqn:rN}) and (\ref{eqn:fixed4}). The solid lines are results of the POS (\ref{eqn:O1}) with $\what\equiv0$ and the symbols are the results of the PS (\ref{eqn:ow}) with $\what\neq0$. The red lines signify the physical positions and phases $r_j^x,g_j$ of the pulses, while the black lines give the differences $r_j^x-r_k^x,r_j^y-r_k^y,g_j-g_k$.}
\label{fig:basin_N4_trajectories}
\end{figure}

In Figure \ref{fig:basin_N45}, we plot the corresponding basins of attraction for (a,b) $N=4$ and (c,d) $N=5$ pulses. The main difference in these cases compared to the $N=3$ scenario, is that now there are different equilibrium solutions possible. For $N=4$ we can identify simulations which integrate to the equilibrium point given by (\ref{eqn:fixed4}), denoted by the black contour, and an example of the time-evolution of the pulse motions for this case is given in Figure \ref{fig:basin_N4_trajectories}(a). In addition to this equilibrium point we also find the green central basin contour which denotes the basin where trajectories integrate to the $N=3$ equilibrium point in (\ref{eqn:fixed3}). In this case, after an initial transient period, the influence of one of the pulses becomes so weak that it remains stationary, while the other three pulses move away as a coherent structure. An example of the time-evolution of this case is given in Figure \ref{fig:basin_N4_trajectories}(b). In this case, after an initial transient period the pulses appear to enter a periodic regime until $\that\approx2\times10^5$ where the $\rb_1,~\rb_2$ and $\rb_3$ pulses begin moving in the negative $x$-direction. The periodic state in Figure \ref{fig:basin_N4_trajectories}(b) is just a transient, but it is also possible to observe persistent periodic limit cycles too, and the basins for these results are denoted by the red contours in Figure \ref{fig:basin_N45}(a,b). These basins tend to be smaller than the other two equilibrium states. This is possibly due to them being more vulnerable to being `knocked' into the other states by small numerical noise, which could be the case for the result plotted in Figure \ref{fig:basin_N4_trajectories}(b). An example of a persisting periodic trajectory can be found in Figure \ref{fig:basin_N4_trajectories}(c). We note that periodic limit cycle solutions do not have to be stationary in physical space, as the result in Figure \ref{fig:basin_N4_trajectories}(c) shows, here the pulses drift slowly in the negative $x$-direction while performing periodic interactions.

The basin contours for both $N=4$ and $N=5$ in Figure \ref{fig:basin_N45} are much more fragmented at the boundaries than the $N=3$ case, suggesting that the system is chaotic close to these points. This is probably not surprising, given the more complex structure of the problems. For $N=5$ in Figure \ref{fig:basin_N45}(c,d) there are now just the two possible basins, which are both fixed equilibrium points, these are the $N=5$ (black contour) and the $N=3$ (green contour) cases in (\ref{eqn:fixed5}) and (\ref{eqn:fixed3}).

For cases with $N>5$, we were not able to identify any stable equilibrium solutions for the parameter ranges we examined with initial configurations similar to those in Figure \ref{fig:initial_fixed_config}, but the results in this section suggest that it is common for larger systems to try to fragment and breakdown into smaller systems, which appear to be more inherently stable. Therefore, in the next section we briefly investigate the evolution of two $N=3$ configurations in the equilibrium point scenario when they are placed next to each other in the plane.

\subsubsection{Interaction of coherent structures}
\label{sec:coherent}

In \S\ref{sec:basins}, we found that as the number of pulses in the system increased, $N>3$, there was an increased possibility that the system would interact in such a way as to `strip off' some pulses and reach an equilibrium state involving only $N=3$ pulses. In Figure \ref{fig:basin_N45}, it was clear that this scenario gave a larger basin of attraction for the initial conditions we used than the respective $N=4$ or $N=5$ basins. In this section, we consider two of these coherent $N=3$ states and investigate their interaction.

\begin{figure}[!htb]
\begin{center}
\includegraphics[width=0.6\textwidth]{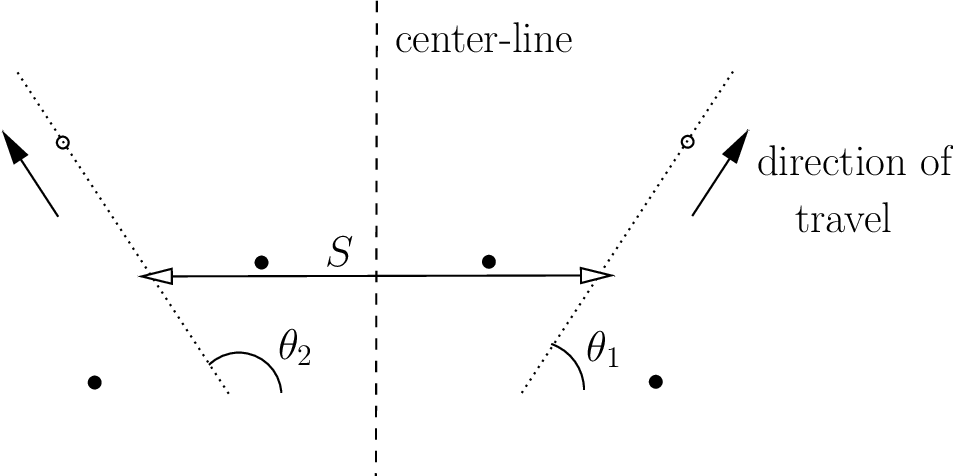}
\end{center}
\caption{Configuration of the two $N=3$ coherent structures from the equilibrium point calculation (\ref{eqn:fixed3}), whose centres are separated by a distance $S$ and each structure is rotated by an angle $\theta_1$ and $\theta_2$.}
\label{fig:2struct_initial}
\end{figure}
The results in Figure \ref{fig:basin_N3_trajectories}(a) and (d) show that the $N=3$ equilibrium point structure propagates in the $(x,y)$-plane in the direction of the pulse which is out-of-phase with the other two, i.e. along the $x$-axis in the example in Figure \ref{fig:fixed_points}(a). With this information, we consider an initial condition of two $N=3$ structures, where the centres of the structures ($\frac{1}{3}(\rb_1+\rb_2+\rb_3)$) are separated by a distance $S$, and each structures' orientation is defined by an angle $\theta_1$ or $\theta_2$, as depicted in Figure \ref{fig:2struct_initial}. These structures will then propagate along the directions indicated in Figure \ref{fig:2struct_initial} and we examine their interaction when/while the pulses of one structure are close enough to affect those in the other structure. In order to simplify the number of parameters in this problem, we only consider the case of symmetrically aligned structures which have reflection symmetry about the centreline, i.e. $\theta_2=\pi-\theta_1$, and we also only consider solutions to the POS in order to identify interesting phenomena, to speed up the calculations. The agreement between the POS and the PS trajectory results in \S\ref{sec:basins} show that this is a robust assumption. For these simulations, we solve the POS using a $4^{\rm th}$ order Runge-Kutta integration scheme, and we use quadruple precision calculations in order to eliminate rounding errors from the long-time runs which can break this symmetry at large times.

We consider the discretised set of initial conditions given by
\[
S_j=4+(j-1)\frac{1}{4},~~~~~~j=[1,21],
\]
and
\[
\theta_{1k}=(k-1)\frac{\pi}{20},~~~~~~k=[1,21],
\]
and run each simulation to $\that=6\times10^4$. We then plot the trajectories and identify by eye what behaviour the structures are exhibiting.

\begin{figure}[!htb]
\begin{center}
\includegraphics[width=0.5\textwidth]{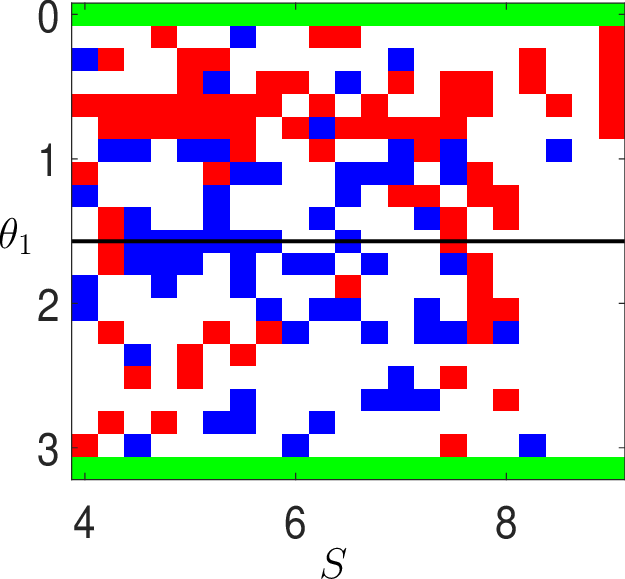}
\end{center}
\caption{Results of the two $N=3$ coherent structure interactions. The red and blue pixels denote results where two 2-pulse pairs are ejected away from the original structure position, the green pixels denote results where the two structures drift slowly towards or away from each other, and the white pixels denotes cases where either the 6 pulses continue to interact with each other, or the two $N=3$ structures remain intact. The red pixels indicate cases where the leading out-of-phase pulse is included in the ejected pair and the blue pixels denotes cases where it is the two in-phase pulses which make up the ejected pair.}
\label{fig:2struct_outcome}
\end{figure}
When  $\theta_1=\pi$ the two structures just propagate away from each other, with a small transient region of adjustment for cases with $S\lesssim5$. For $\theta_1=0$ the two structures propagate towards each other, reaching a slow drifting state when the two leading pulses are a distance of $\approx2.6$ units apart. At this point the two trailing pulses slowly drift outward and towards the centreline, while for the other combinations of $\theta_1$ and $S$, we observe scenarios where the pulses either remain as 6 interacting pulses, two sets of 3 interacting pulses, or we observe two-sets of 2 pulses being stripped away to propagate away from the initial system. This final case is another example of a large system of pulses wanting to degenerate into smaller systems.
\begin{figure}[!htb]
\begin{center}
(a)\includegraphics[width=0.4\textwidth]{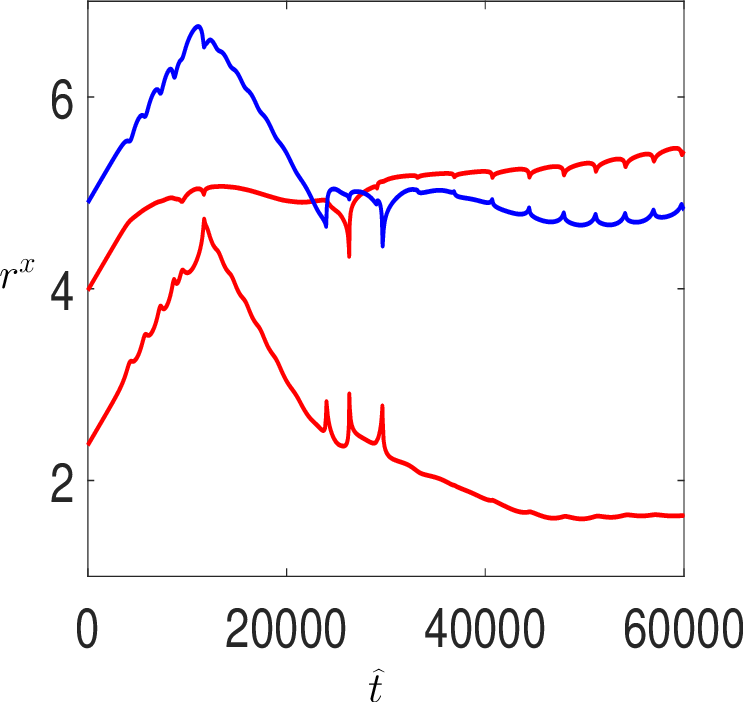}
\includegraphics[width=0.4\textwidth]{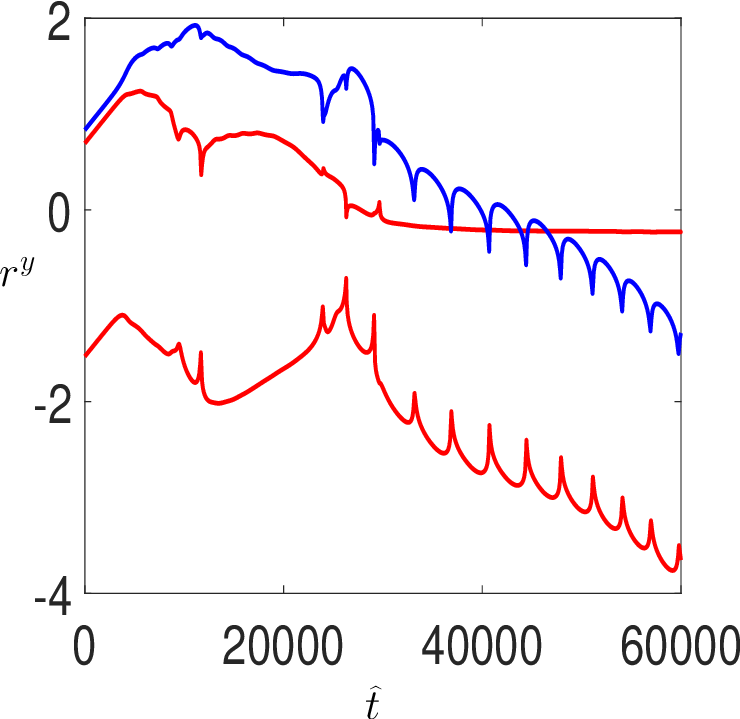}
(b)\includegraphics[width=0.4\textwidth]{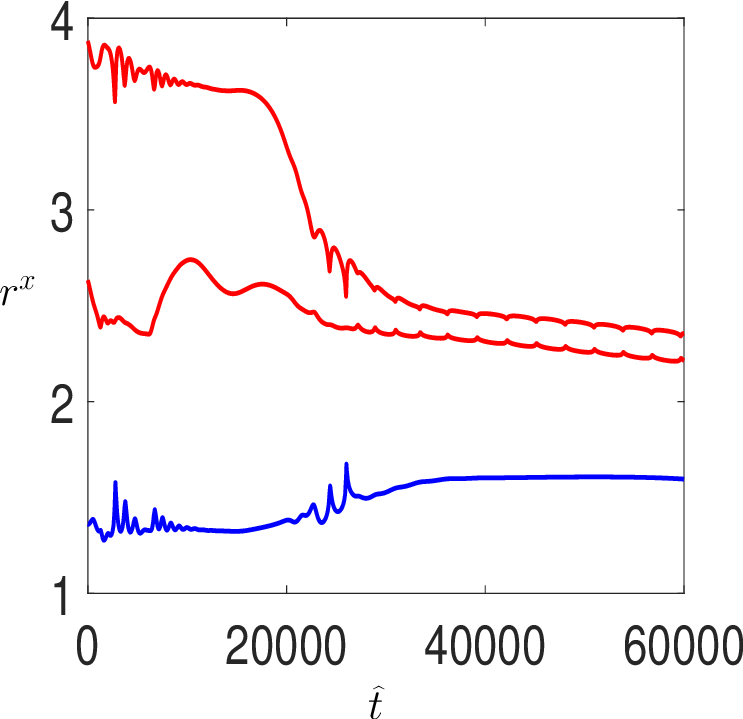}
\includegraphics[width=0.4\textwidth]{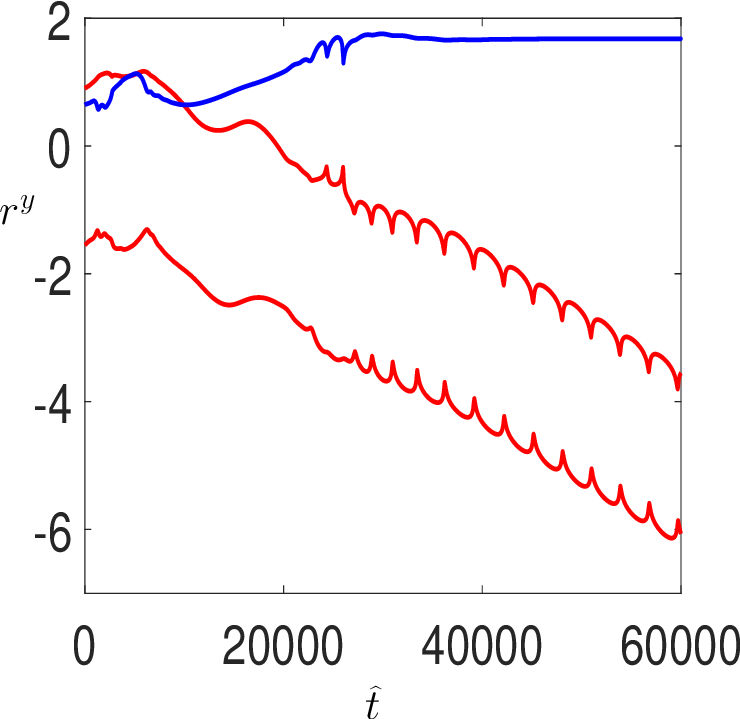}
\end{center}
\caption{Time-evolution plots of the pulse parameters from the POS for the rightmost set of pulses from Figure \ref{fig:2struct_initial} for the initial conditions (a) $S=7.5$, $\theta_1=4\pi/5$ and (b) $S=5.25$, $\theta_1=3\pi/20$. In each case the red lines are the two rear pulses which are initially in phase with each other, while the blue line is the trajectory of the initially out-of-phase pulse.}
\label{fig:2struct_trajectories}
\end{figure}
In Figure \ref{fig:2struct_outcome}, we identify where in the $(S,\theta_1)$-parameter space these pairs of propagating pulses exist. In Figure \ref{fig:2struct_outcome}, the red and blue pixels indicate parameter values where the  pair of pulses is ejected from the initial configuration, after some initial  transient period. Two examples of this can be seen in the pulse position plots in Figure \ref{fig:2struct_trajectories}. The difference between the colour of the pixels is, the red pixel denotes cases where the leading pulse is included in the pair of pulses propagating away (see Figure \ref{fig:2struct_trajectories}(a)), while for the blue pixels the pair of pulses propagating away are the rear two pulses (see Figure \ref{fig:2struct_trajectories}(b)). What Figure \ref{fig:2struct_outcome} shows is that there is no clear, coherent region of the $(S,\theta_1)$-plane where these two pulses are ejected from the system, and it is strongly dependent on the initial condition. Clearly, when the two structures are initially closer, and moving towards each other ($\theta_1\leq\pi/2$) there are more interactions which ejects the two pulse pairs. However, even for the case $S=7.5$ with $\theta_1=4\pi/5$, such that the structures are moving apart, as shown in Figure \ref{fig:2struct_trajectories}(a), two of the pulses are ejected. It is plausible that in those cases where two pulses are not ejected (white pixels in Figure \ref{fig:2struct_outcome}) they might eventually lead to this scenario if the problem is integrated for longer. However, in these cases, numerical round-off becomes significant and can contaminate the results.

\section{Conclusions and Discussion}
\label{sec:conclusions}

In this paper, we examined the dynamics of multi-pulse interactions in two-space dimensions in the quintic-complex Ginzburg-Landau equation (QCGLE). The approach used was based on the centre-manifold approach used on the one-dimensional equivalent problem in \cite{rossides2023}. In this approach, the governing equations were formulated  by writing the solution as a sum of the individual pulses plus a remainder function. The dynamics of the pulses are then determined by projecting the solution onto the stable centre-manifold such that a system of {\it slow} ODEs govern the pulse motion, coupled to a {\it fast} PDE for the remainder function. This Projected Scheme (PS) is advantageous over more direct numerical methods to solve the QCGLE because it has a distinct speed advantage, which is only restricted by the time taken to solve the stationary approximation of the fast PDE. In this work, we also highlight an asymptotic approach which allows us to investigate these pulse interactions when these pulses are very-well-separated, and in this case the governing PDE for the remainder function reduces to a static PDE, which can be solved readily using GMRES methods.

Initially, we applied the centre-manifold reduction scheme to the case of two interacting pulses. In this case, the pulses interact along a line joining the centre of the pulses, so the problem reduces to a quasi-one-dimensional problem. We showed that the problem shares the same solution structure as the one-dimension problem in \cite{rossides2023}. Namely, that when neglecting the remainder function, $w(x,y)$, the so called Projected ODE Scheme (POS) limit, the phase plane $(\rbar\cos\gbar,\rbar\sin\gbar)$ contains cells of periodic orbits circling a centre. When the remainder function  is included for a typical parameter set with $\beta_r={\rm Re}(\beta)=-0.02$ then the closed periodic orbits open up and the centre becomes a stable focus with all the trajectories spiraling into this point. As $\beta_r$ is varied, the form of the underlying pulse solution changes and for the case when $\beta_r=-2$, we were able to identify a stable limit cycle in this cell to which all initial trajectories converged to.

We then considered the case of $N\geq3$ pulses in the two-dimensional plane. For each of the cases $3\leq N\leq5$ we were able to identify stable equilibrium solutions where the pulses were arranged in near regular geometric configurations. For $N=3$ we find a large basin of attraction for this equilibrium point using the POS and we confirm its robustness via the PS. This shows some differences in the trajectories near to the basin boundary where the pulses are initially close together, such that the assumptions in the POS are invalidated. For $N=4$, we find a much smaller basin of attraction for the 4 pulse equilibrium point, and in fact we find a larger basin for initial conditions which converge to the $N=3$ pulse equilibrium point. In these cases, one of the initial pulses is ejected from the system leaving just three interacting pulses. In the $N=4$ case we also identified solutions which converged to periodic limit cycles, as we found in the $N=2$ case. For $N=5$ pulses, we find just the $N=3$ and $N=5$ equilibrium points, without any periodic limit cycles.

Finally, we investigated the interaction of two $N=3$ stable coherent structures. Here each coherent structure moved in the $(x,y)$-plane in the direction of the pulse which is out-of-phase with the other two. We found that when these two structures are symmetrically placed together they can degenerate into systems which have two sets of two ejected pulses, which propagate away from the initial configuration, showing that the QCGLE multi-pulse systems can be relatively unstable to high-order coherent structures.

The motivation for considering the interaction of coherent structures in this project was the observation that coherent hexagonal structures are observed in the ferrofluid experiment \cite{lloyd2015} and also in the Swift-Hohenberg equation \cite{lloyd2021} with a quadratic and cubic nonlinearity, which is a canonical equation capturing many fundamental features of the ferrofluid problem. In both these cases the hexagonal patterns of pulses are stationary and hence it would be easier to examine their interaction as they would not be constantly propagating and affecting their overall dynamics. These hexagonal patches also appear to be stable, as they are observable  in the experiments, and so perturbing them with another hexagonal patch might lead to further stable structures, including the possibility of them coalescing. Hence, one interesting avenue of future research would be to apply the centre-manifold reduction technique on the Swift-Hohenberg equation and examine pulse interactions in this context.

\red{The centre-manifold projection method in this paper can be readily extended to any dissipative partial differential equation of the form (\ref{eqn:disppde}) of arbitrary spatial dimension, $P$. In the case of symmetric $P$-dimensional pulses, the steady problem remains a function of $r=(x_1^2+x_2^2+\cdots+x_P^2)^{1/2}$ only, and is tractable, along with the eigenvalue problem for the $P+1$ independent eigenmodes. The resulting projection onto the eigenmodes leads to a system of ODEs of size $PN$ for the $PN$ unknown system symmetries $\boldsymbol\xi(t)$. Finally, the remainder function $w(\mathbf{x},t)$ is a solution to a $P$-dimensional linear partial differential equation, and while in the weakly-interacting limit this can be written in the form (\ref{eqn:ow}), the fast and effective solution for $w$ is the limiting factor on the generalization of this approach.}

\vspace{0.25cm}

\section*{Acknowledgements}

This research was funded in part, by the Engineering and Physical Sciences Research Council via grant number UKRI070. For the purpose of Open Access, the authors have applied a Creative Commons Attribution (CC BY) public copyright licence to any Author Accepted Manuscript version arising from this submission.

\vspace{0.25cm}

\section*{Data Availability Statement}

The main code used to generate the results in figures 4,5,6 can be found on Github at \url{https://github.com/David-JB-Lloyd-Lab/Multi_Pulse_Complex_Ginzburg_Landau_Equation_2D}

\appendix

\section{Statement of Theorem 2.1 from \cite{rossides2023}}
\label{appen:theorem}

This theorem was first stated and proved in \cite{ZelMiel09}.

Assume the pulse solution $V$ satisfies the QCGLE (\ref{eqn:gl}), and the PDEs linear operator $\mathbb{L}$, given in (\ref{eqn:lineargl}), is such that it has distinct zero eigenvalues (no Jordan blocks are allowed at zero), and the rest of the spectrum of $\mathbb{L}$ lies in the stable half plane. Then for all separation distances $\widehat{d}\geq d$, for sufficiently large $d$, there exists a small neighbourhood of
\[
\mathscr{P}_n=\mathscr{P}_n(\widehat{d}):=\left\{V_{\bX}({\bf x}), \bX\in(\mathbb{R}\times\mathbb{R}\times\mathbb{S}_1), \widehat{d}\geq d\right\},
\]
where
\[
V_{\bX}({\bf x})=\sum_{k=1}^N V_k({\bf x}),~~~~~{\bf x}=(x,y),
\]
(which is independent of $\widehat{d}$), say,  in the metric of $L^{\infty}(\mathbb{R})$ such that the projected system \ref{eqn:Ow} together with \ref{eqn:w} is equivalent to the initial QCGLE (\ref{eqn:gl}) as long as the trajectory $u({\bf x},t)$ remains in this neighbourhood.

Moreover, there exists a $C^2$-smooth function $W:\overline{\mathscr{P}_n(\widehat{d})}\to L^\infty(\mathbb{R})$ such that
\[
||W||_{C^2}\leq C\epsilon,~~~~\epsilon:=e^{-\lambda_r d},
\]
where $\lambda_r$ is from (\ref{eqn:disp}), and the manifold defined by the graph of $W$ ($w=W(\bX),~\bX\in\mathscr{P}_n(\widehat{d})$) is an invariant centre manifold of the projected system. In particular, the dynamics on this manifold is determined by the system of ODEs (\ref{eqn:matrix2}) with 
\[
||{\bf C}^{-1}(\epsilon,\bX,W(\bX)){\bf F}(\epsilon,\bX,W(\bX))||_{C^2}\leq C\epsilon,
\]
and the trajectory $\bX$ starting on this manifold can only leave through the boundary $\partial\mathscr{P}_n(\widehat{d})$.

Finally, this manifold is exponentially stable and normally hyperbolic, so for any trajectory $(\bX(t),w(t))$ of the projected system there is a trace trajectory $(\bX_0(t),W(\bX_0(t)))$ on the manifold such that
\[
||\bX(t)-\bX_0(t)||+||w(t)-W(\bX_0(t))||_{L^\infty}\leq C e^{-\kappa t},
\]
for some positive $\kappa$ which is independent of $d$. This estimate holds until the trace trajectory $\bX_0$ reaches the boundary $\partial\mathscr{P}_n(\widehat{d})$.


\bibliographystyle{plain}
\bibliography{CGLpaperbib.bib}


\end{document}